\documentclass[preprint,11pt]{elsarticle}
\usepackage{mathrsfs}
\usepackage{amsfonts}
\usepackage{amsmath}
\usepackage{amstext}
\usepackage{stmaryrd}
\usepackage{amssymb}
\usepackage{amsthm}
\usepackage{mathrsfs}
\usepackage{amsfonts}
\usepackage{enumitem}
\usepackage{amscd}
\usepackage{indentfirst}
\usepackage{amsmath,amsfonts,amssymb,amsthm}
\usepackage{amsmath,amssymb,amsthm,amscd}
\usepackage{graphicx,mathrsfs}
\usepackage{subcaption}
\usepackage{appendix}
\usepackage[colorlinks,urlcolor=red]{hyperref}
\usepackage[normalem]{ulem}

\newtheorem{theorem}{Theorem}[section]
\newtheorem{lemma}[theorem]{Lemma}
\newtheorem{definition}[theorem]{Definition}
\newtheorem{proposition}[theorem]{Proposition}
\newtheorem{assum}{Assumption}

\newcommand{\subscript}[2]{$#1 _ #2$}
\newcommand{\mlb}{\mathcal{L}^B_{\kappa}}
\newcommand{\veps}{\varepsilon}
\newcommand{\dlt}{\delta}
\newcommand{\ptl}{\partial}
\newcommand{\bx}{\boldsymbol{x}}%bblack
\newcommand{\gx}{\boldsymbol{X}}%ggreat

\newcommand{\bk}{\mathbf{k}}
\newcommand{\gk}{\mathbf{K}}
\newcommand{\frk}{\mathfrak{K}}
\newcommand{\bm}{\mathbf{m}}
\newcommand{\bn}{\mathbf{n}}

\newcommand{\bv}{\mathbf{v}}
\newcommand{\br}{\mathbb{R}}
\newcommand{\balpha}{\boldsymbol{\alpha}}
\newcommand{\bbeta}{\boldsymbol{\beta}}
\newcommand{\ml}{\mathcal{L}}
\newcommand{\bla}{\big\langle}
\newcommand{\bra}{\big\rangle}
\newcommand{\bav}{\Big|}
\newcommand{\bkappa}{\boldsymbol{\kappa}}

\numberwithin{equation}{section}

\theoremstyle{definition}

\begin{document}
\begin{frontmatter}

\title{Wave packet dynamics in slowly modulated photonic graphene
%\tnoteref{sup}} \tnotetext[sup]{This work is partially supported by the NSFC under the grants 11371282
}
\author[cam]{Peng Xie}
\ead{xiep14@mails.tsinghua.edu.cn}

\author[cam]{Yi Zhu\corref{cor1}}
\ead{yizhu@tsinghua.edu.cn}
\cortext[cor1]{Corresponding author}

%\ead{email address}
%\fntext[label2]{}
%\address{Address\fnref{label3}}
%\fntext[label3]{}

%\author[cam,csrc]{Peng Xie}
%\ead{xiep14@mails.tsinghua.edu.cn}
%\author[cam]{Yi Zhu}

\address[cam]{Zhou Pei-Yuan Center for Applied Mathematics, Tsinghua University, Beijing 100084, China}


\begin{abstract}
Mathematical analysis on electromagnetic waves in photonic graphene, a photonic topological material which has a honeycomb structure, is one of the most important current research topics. By modulating the honeycomb structure, numerous
topological phenomena have been observed recently.  The electromagnetic waves in such a media are generally described by the 2-dimensional wave equation. It has been shown that the corresponding elliptic operator
with a honeycomb material weight has Dirac points in its dispersion surfaces. In this paper, we study the time evolution of the wave packets spectrally concentrated at such Dirac points in a modulated honeycomb material weight.
We prove that such wave packet dynamics is governed by the Dirac equation with a varying mass in a large but finite time. Our analysis provides mathematical insights to those topological phenomena in photonic graphene.
\end{abstract}
\end{frontmatter}

%\addcontentsline{toc}{section}{References}
%\newpage

%\tableofcontents
%\section{Introduction}
%\addcontentsline{toc}{section}{Introduction}
\section{Introduction}
This paper is concerned with a two dimensional (2D) wave equation
\begin{equation} \label{weq:perturbed}
	\partial_{t}^2\psi-\nabla \cdot \big( W_{\veps}(\bx)\nabla \psi\big) = 0,\quad (\bx\in\br^2,~t>0)
\end{equation}
where the material weight $ W_{\veps}$ has a slowly-modulated honeycomb structure,
\begin{equation}
W_{\veps}(\bx)=A(\bx)+\veps\kappa(\veps\bx)B(\bx).
\end{equation}
Here $A(\bx)$ defines a honeycomb structured media,  $\veps$ is a sufficiently small positive parameter,  $\kappa(\veps\bx)$ together with $B(\bx)$ are modulations applied to the honeycomb material, (see Section 2 for details). The 2D wave equation arises in describing the electromagnetic waves in an optical metamaterial whose permittivity and/or permeability  are designed to have a hexagonal symmetry.  These materials are also referred to as ``photonic graphene". In applications, modulations are often added to the honeycomb structure, for instance, domain-wall-like defects. These modulations leads to numerous topological phenomena \cite{drouot2018defect,haldane2008Possible,khanikaev2013photonic,lee2018elliptic,lu2018dirac,lu2014Topological}

The subtle property of a honeycomb structured media defined by $A(\bx)$ is the existence of two-fold degenerate points, called Dirac points $(\gk_*,E_D)$, lying in the dispersion surfaces of the elliptic operator $\mathcal{L}^A\equiv -\nabla \cdot \big(A(\bx)\nabla\big)$. The corresponding eigenfunctions, $\Phi_1(\bx)$ and $\Phi_2(\bx)=\overline{\Phi_1(-\bx)}$ satisfy
\begin{equation}
\ml^A \Phi_j(\bx)=E_D\Phi_j(\bx),\quad j=1,~2.
\end{equation}
Due to the conical property of Dirac points, the wave packets associated with such points have very novel propagation patterns. The goal of this work is to investigate the effective dynamics of such wave packets in \eqref{weq:perturbed} by deriving the envelope equation with a rigorous justification. To this end, we consider the wave equation \eqref{weq:perturbed} with the initial condition
\begin{equation}\label{generalinitial}
\left\{
\begin{split}	
\psi(\bx,0) =&~ \veps\Big[a_1(\veps \bx)\Phi_1(\bx)+a_2(\veps \bx)\Phi_2(\bx)\Big], \\
\ptl_t\psi(\bx,0) =&~ \veps \Big[b_1(\veps \bx)\Phi_1(\bx)+b_2(\veps \bx)\Phi_2(\bx)\Big],
\end{split}
\right.
\end{equation}
where $a_j(\cdot),~b_j(\cdot),~j=1,~2$ are all in Schwartz space $\mathcal{S}(\br^2)$. This initial condition represents a general wave-packet spectrally localized at Dirac point $(\gk_*,E_D)$. Thanks to the linearity of \eqref{weq:perturbed}, the solution has two branches which evolve independently. Thus it is sufficient to consider the initial condition corresponding to one of the branches \footnote{The analysis for the other branch is essentially the same by just changing $\sqrt{E_D}$ to $-\sqrt{E_D}$.}
\begin{equation}\label{ini}
\left\{
\begin{aligned}
	\psi(\bx,0) =&~ \veps\Big[\alpha_{10}(\veps \bx)\Phi_1(\bx)+\alpha_{20}(\veps \bx)\Phi_2(\bx)\Big], \\
	\ptl_t\psi(\bx,0) =&~ i\sqrt{E_D}\veps \Big[\alpha_{10}(\veps \bx)\Phi_1(\bx)+\alpha_{20}(\veps \bx)\Phi_2(\bx)\Big],
\end{aligned}
\right.
\end{equation}
where $\alpha_{10}(\cdot),~\alpha_{20}(\cdot)$  are in Schwartz space $\mathcal{S}(\br^2)$. We remark that the factor $\veps$ in front of $\psi(\bx,0)$ and $\ptl_t \psi(\bx,0)$ are not essential due to the linearity of our problem. But this choice ensures that which brings great convenience to our analysis, i.e., $\|\psi(\bx,0)\|_{H^s(\br^2)}=\mathcal{O}(1)$, $\|\ptl_t\psi(\bx,0)\|_{H^s(\br^2)}=\mathcal{O}(1)$.

We shall prove that the wave equation \eqref{weq:perturbed} with initial condition \eqref{ini} has the following asymptotic solution,
\begin{eqnarray}\label{solu1}
    \psi(\bx, t) = e^{i\sqrt{E_D}t}\veps\Big[\alpha_{1}(\veps \bx,\veps t)\Phi_1(\bx)+\alpha_{2}(\veps \bx,\veps t)\Phi_2(\bx)\Big]+\eta(\bx, t),
\end{eqnarray}
where $\alpha_j(\veps\bx, \veps t),~j=1,~2$ represent the slowly varying envelopes and $\eta(\bx, t)$ is the small correction to the leading order approximation. The main result of this paper stated in Theorem \ref{mainthm} is as follows. The envelope $\alpha_j(\gx,T),~ j=1,~2~\big(\gx=(X_1,~X_2)=\veps\bx,~T=\veps t\big)$ satisfy the so-called Dirac equation with a varying mass
\begin{equation}\label{dirac}
	\left\{
	\begin{aligned}
		i\ptl_{T}\alpha_1-\frac{v_{_F}}{2\sqrt{E_D}}\big(i\ptl_{X_1}-\ptl_{X_2}\big)\alpha_2+\frac{\vartheta_{\sharp}\kappa(\gx)}{2\sqrt{E_D}}\alpha_1=0 \\
		i\ptl_{T}\alpha_2-\frac{v_{_F}}{2\sqrt{E_D}}\big(i\ptl_{X_1}+\ptl_{X_2}\big)\alpha_1-\frac{\vartheta_{\sharp}\kappa(\gx)}{2\sqrt{E_D}}\alpha_2=0
	\end{aligned}~,
	\right.
\end{equation}
with the initial condition
\begin{equation}\label{dirac:ini}
	\alpha_{1}(\gx, 0)=\alpha_{10}(\gx),~\alpha_{2}(\gx, 0)=\alpha_{20}(\gx).
\end{equation}
Furthermore, given $\veps>0$ sufficiently small, $\eta(\bx, t)$ can be controlled in Sobolev space $H^s(\br^2)$ over a large but finite time scale $t\in[0,~\rho\veps^{^{-1}}]$ as follows
\begin{equation}\label{conclusion}
\sup_{0\le t \le \rho\veps^{^{-1}}} \|\eta(\bx, t)\|_{H^s(\br^2)}\le C\veps^{1-\nu},
\end{equation}
for any $0<\nu<1,~s\ge0$ and $\rho>0$.

The main idea of our proof is inspired from the pioneer work by Fefferman and Weinstein \cite{fefferman2014wave}. They gave a rigorous justification of massless Dirac equation derived from the Schr\"odinger equation with a perfect honeycomb potential  which corresponds to the special case  $\kappa(\veps \bx)\equiv0$. In our current work, the case where $\kappa(\veps \bx)\neq 0$ is considered in the wave equation \eqref{weq:perturbed} which has a second order derivative in time. These differences will bring technical difficulties in analysis which require new treatments.

In past decades of years, one of the most popular research subject is to understand and realize topological phenomena in different topological materials. One of the most successful example is the honeycomb-based material \cite{geim2007The, neto2009electronic, rechtsman2013Photonic}. It stimulates the mathematical analysis of the Schr\"{o}dinger equation with a honeycomb potential over the past few years \cite{ablowitz2012tight, ablowitz2009Conical, arbunich2018rigorous, birman2006homogenization, fefferman2012honeycomb, fefferman2014wave}. Fefferman and Weinstein rigorously proved the existence of Dirac points of Schr\"odinger operator with a generic honeycomb potential \cite{fefferman2012honeycomb}, and later gave a mathematical justification of the massless Dirac equation which governs the dynamics of the wave packets associated with Dirac points \cite{fefferman2014wave}. Topological edges states, strong binding limit, nonlinearity, and other aspects regarding this equation have also been investigated by them and others \cite{ablowitz2012nonlinear, ablowitz2013Nonlinear2, arbunich2018rigorous, fefferman2016edge, fefferman2017topologically, fefferman2018honeycomb, fefferman2018edge, giannoulis2008interaction, keller2018spectral}. Meanwhile, many ``artificial graphenes'', analogies of graphene in other fields, have been created to realize similar properties. Amongst those, photonic graphene has attracted a lot of interest due to its potential applications and relatively simple experimental realizations \cite{haldane2008Possible, khanikaev2013photonic, lu2014Topological, raghu2008analogs}. To study electromagnetic waves in photonic graphene, we have to deal with Maxwell's equations in \cite{de2014effective, de2017derivation}. In a simple physical setting, for example, the propagation of transverse electronic fields can be described by the 2D wave equation \eqref{weq:perturbed}. In a previous work, Lee-Thorp, Weinstein and Zhu \cite{lee2018elliptic} rigorously proved the existence of Dirac points of 2D elliptic operator with a honeycomb structured material weight and the existence of topological edge states in a domain-wall-modulated honeycomb structure. Their work paved the way to the mathematical analysis of electromagnetic waves in topological photonic materials. However, the wave packet dynamics has not been thoughtfully studied yet. In this present work, we shall give a rigorous investigation along this interesting and important line.

Before proceeding, we use the following basic requirements on the material weight.
\begin{definition}\label{assume}
The $2\times2$ complex-valued matrix function $W(\bx)$ is called a material weight, if
\begin{enumerate}
	\item $W(\bx)$ is Hermitian and smooth for all $\bx\in\br^2$,
	\item  $W(\bx)$ is elliptic, i.e., for any $\boldsymbol{\xi}\in \mathbb{C}^2$, $\exists~0<C_1\le C_2<+\infty$, such that $C_1|\boldsymbol{\xi}|^2\le\bar{\boldsymbol{\xi}}\cdot W(\bx)\boldsymbol{\xi}\le C_2|\boldsymbol{\xi}|^2$ for all $\bx\in\br^2$.
\end{enumerate}
\end{definition}

The rest of the paper is organized as follows. In section 2, we briefly review the basic Floquet-Bloch theory for 2D periodic elliptic operator, honeycomb structured media, Dirac points and other results which are used in later analysis; In section 3, we present the well-posedness of the envelope equation \eqref{dirac} in Schwartz space and its connection to topological edge states; In section 4, we conclude the main result$-$Theorem \ref{mainthm}; In section 5, we give the detailed proofs of the key estimates which are essential to the proof of Theorem \ref{mainthm}. In the appendix, we discuss how to apply our analysis to the non-modulating case, i.e., $\kappa(\veps\bx)\equiv0$, and give the detailed estimate of the solution to Dirac equation \eqref{dirac} in Schwartz space.
\newline

The following notations and conventions are used in this work:
\begin{enumerate}
	\item $\Omega$ is the fundamental cell in \eqref{cell}, and $\Omega^*$ the dual fundamental cell defined in \eqref{dualcell}.
	\item The operators $\ml^A,~\ml^B,~\mlb$ are denoted as follows:$$\ml^A=-\nabla\cdot\big(A(\bx)\nabla\big),~\ml^B=-\nabla\cdot\big(B(\bx)\nabla\big),~\mlb=-\nabla\cdot\big(\kappa(\veps\bx)B(\bx)\nabla\big).$$
	\item The standard notations for function spaces are used:
	\begin{itemize}
		\item[(1)] $H^s(\br^2)$ is the Sobolev space, i.e., $f(\bx)\in H^s(\br^2)$, $\|f\|_{H^s(\br^2)}^2=\sum\limits_{|\bn|\le s}\int_{\br^2}|\ptl^{\bn}_{\bx}f(\bx)|^2d\bx<\infty$;
		\item[(2)] $\mathcal{S}(\br^2)$ is the Schwartz space, i.e., $f(\bx)\in\mathcal{S}(\br^2)$,  $\sup\limits_{\bx\in\br^2}|\bx^{\bm}\ptl_{\bx}^{\bn}f(\bx)|<\infty,~\forall~\bm,\bn\in\mathbb{N}^2$;
		\item[(3)] $C_b^{\infty}(\br^2)$ contains all functions which are bounded and smooth.		
	\end{itemize}
	\item We use the indicator function $\chi(D)$ to define the Bloch spectral cutoff, i.e., for any constant $C\ge0$, and $f:\Omega^*\to\mathbb{C}$,
		\begin{equation*}
			\chi(|\bk|\le C)f(\bk):=\left\{
				\begin{aligned}
					f(\bk),\quad\text{if}~|\bk|\le C, \\
					0,\quad\text{otherwise}.
				\end{aligned}
				\right.
		\end{equation*}
	\item The Pauli matrices $\sigma_j,~j=1,~2,~3$ are:
		\begin{equation*}
		\sigma_1=
		\begin{pmatrix}
		0 & 1 \\
		1 & 0
		\end{pmatrix},\quad
		\sigma_2=
		\begin{pmatrix}
		0 & -i \\
		i & 0
		\end{pmatrix},\quad
		\sigma_3=
		\begin{pmatrix}
		1 & 0 \\
		0 & -1
		\end{pmatrix}.
		\end{equation*}	
	\item $A\approx B$ if and only if there exist two constants $0<C_1\le C_2$ such that $C_1 A\le B \le C_2 A$.
	
	\item We use notations $ \bla \cdot, \cdot\bra$ and $\bla \cdot, \cdot\bra_{L^2(\Omega)}$ to distinguish the inner products on $L^2(\br^2)$ and $L^2(\Omega)$, i.e.,
		\begin{equation*}
			\bla f_1(\bx), f_2(\bx)\bra=\int_{\br^2}\overline{f_1(\bx)}f_2(\bx) d\bx,\qquad
		\bla g_1(\bx), g_2(\bx)\bra_{L^2(\Omega)}=\int_{\Omega}\overline{g_1(\bx)}g_2(\bx) d\bx.
		\end{equation*}
	\item For any $f(\bx)\in L^2(\br^2)$, $\hat{f}(\boldsymbol{\xi})$ represents its Fourier transform while $\tilde{f}_b(\bk)$ stands for its Bloch component, i.e.,
		\begin{equation*}
			\hat{f}(\boldsymbol{\xi})=\int_{\br^2}e^{-i\boldsymbol{\xi}\cdot\bx}f(\bx) d\bx,\qquad \tilde{f}_b(\bk) = \bla\Phi_{b}(\cdot;\bk),~f(\cdot)\bra~.
		\end{equation*}
     \item The repeated index summation convention is used throughout.
\end{enumerate}

\section{Preliminaries}

In this section, we list the Floquet-Bloch theory, honeycomb structure, Dirac points which are desired for the arguments of this work. We refer readers to \cite{bensoussan1978asymptotic, bloch1929u, conca1997homogenization,eastham1973spectral, kuchment2012floquet, lee2018elliptic, pelinovsky2011localization, wilcox1978theory} for more details.

\subsection{Floquet-Bloch theory}

A lattice $\Lambda$ in $\br^2$ is generated by two linearly independent vectors $\bv_1$, $\bv_2$, i.e.,
\begin{eqnarray}
\Lambda = \{\bm\bv=m_1\bv_1+m_2\bv_2: \bm=(m_1, m_2)\in\mathbb{Z}^2\}=\mathbb{Z}\bv_1\oplus\mathbb{Z}\bv_2,
\end{eqnarray}
and the fundamental cell is chosen to be the parallelogram:
\begin{eqnarray}\label{cell}
\Omega = \{\theta_1\bv_1+\theta_2\bv_2: 0\le\theta_j\le 1,~ j=1,~ 2\},
\end{eqnarray}
where $|\Omega|$  is the area of $\Omega$. In this work, we specify the lattice $\Lambda$ to be a triangular lattice with lattice vectors
\begin{equation}\label{latticevectors}
\bv_1=
\begin{pmatrix}
\frac{\sqrt{3}}{2}\\
\frac12
\end{pmatrix}
,~~~\bv_2=
\begin{pmatrix}
\frac{\sqrt{3}}{2}\\
-\frac12
\end{pmatrix}.
\end{equation}

The dual lattice
\begin{equation}
	\Lambda^* = \{\bm\bk=m_1\bk_1+m_2\bk_2: \bm=(m_1, m_2)\in\mathbb{Z}^2\}=\mathbb{Z}\bk_1\oplus\mathbb{Z}\bk_2~,
\end{equation}
is generated by the dual lattice vectors $\bk_1,~\bk_2$ which satisfy $\bk_i\cdot\bv_j = 2\pi\delta_{ij},~i,~j=1,~2$. For the triangular lattice defined by \eqref{latticevectors}, the dual lattice vectors are
\begin{equation}
  \bk_1=\frac{4\sqrt{3}}{3}\pi
    \begin{pmatrix}
      \frac12\\
      \frac{\sqrt{3}}{2}
    \end{pmatrix}
    ,\quad\bk_2=\frac{4\sqrt{3}}{3}\pi
    \begin{pmatrix}
      \frac12\\
      -\frac{\sqrt{3}}{2}
    \end{pmatrix}.
\end{equation}
Throughout this work, we choose the parallelogram $\Omega^*$:
\begin{equation}\label{dualcell}
\Omega^* = \{\theta_1\bk_1+\theta_2\bk_2: -\frac12\le\theta_j\le \frac12,~ j=1,~2\},
\end{equation}
as the fundamental dual cell.\footnote{In many literatures, another choice of the fundamental cell is the Brillouin zone $\mathcal{B}$, consisting the points $\bk \in \br^2$ which are closer to the origin than to any other lattice points in $\Lambda^*$. For the triangular lattice, the Brillouin zone is a hexagon, see for example \cite{fefferman2012honeycomb}. }
\newline

Let $L^2_{\text{per}}(\br^2/\Lambda)$ denote a subspace of  $L^2_{\text{loc}}(\br^2)$ containing all $\Lambda$-periodic functions, namely $f\in L^2_{\text{loc}}(\br^2)$ and $f(\bx+\bv)=f(\bx)$, $\forall~\bx\in\br^2,~\forall~\bv\in\Lambda$. For each $\bk\in\br^2$, we denote $f\in L_{\bk}^2(\br^2/\Lambda)$ if $e^{-i\bk\cdot\bx}f(\bx)\in L^2_{\text{per}}(\br^2/\Lambda),~i.e.,~f(\bx+\bv)=e^{i\bk\cdot\bv}f(\bx)$, for all $\bx\in\br^2,~\bv\in \Lambda$. Similarly, we can define the Sobolev space $H^s_\bk(\br^2/\Lambda),~ s\ge0$.

Suppose the matrix function $A(\bx)$ is a material weight in the sense of Definition \ref{assume}, and further $\Lambda$-periodic. $\ml^A=-\nabla\cdot \big(A(\bx)\nabla\big)$ is an elliptic operator with periodic coefficient, thus the Floquet-Bloch theory applies. For any $\bk\in \br^2$, consider the following $L^2_{\bk}$-Floquet-Bloch elliptic eigenvalue problem with pseudo-periodic boundary condition
\begin{align}
		\mathcal{L}^A\Phi(\bx;\bk) &= E(\bk)\Phi(\bx;\bk),\quad\forall~\bx\in\br^2, \\
		\Phi(\bx+\bv;\bk) &= e^{i\bk\cdot\bv}\Phi(\bx;\bk),\quad\forall~ \bv\in\Lambda.
\end{align}
Since the above eigenvalue problem is invariant under the translation $\bk\rightarrow\bk+\bk'$, $\forall~\bk'\in\Lambda^*$, it is sufficient to just pay attention to $\bk$ varying over $\Omega^*$. Alternatively, one can obtain the periodic elliptic boundary problem by setting $\Phi(\bx;\bk)=e^{i\bk\cdot\bx}\phi(\bx;\bk),~\bk\in \Omega^*$, i.e.,
\begin{align}
	\mathcal{L}^A(\bk)\phi(\bx;\bk) &= E(\bk)\phi(\bx;\bk),\quad\forall~\bx\in \br^2, \\
	\phi(\bx+\bv;\bk) &= \phi(\bx;\bk),\quad\forall~\bv\in\Lambda,
\end{align}
where
$$\mathcal{L}^A(\bk)=e^{-i\bk\cdot\bx}\mathcal{L}^A e^{i\bk\cdot \bx}=-(\nabla+i\bk)\cdot \big(A(\bx)(\nabla+i\bk)\big).$$

Notice that $\ml^A(\bk)$ is a self-adjoint elliptic operator in $L^2_{\text{per}}(\br^2/\Lambda)$.  For each fixed $\bk\in \Omega^*$, the above eigenvalue problem has a series of discrete spectrum (or eigenvalues) \cite{bensoussan1978asymptotic,fefferman2014wave,wilcox1978theory}:
\begin{equation}\label{sequence:eigenva}
	0\le E_1(\bk)\le E_2(\bk)\le E_3(\bk)\le\cdots,
\end{equation}
with eigenpairs $\big(\phi_b(\bx;\bk),~E_b(\bk)\big),~b\ge1$, where $\{\phi_b(\bx;\bk)\}_{b\ge1}$ can be chosen a complete orthogonal basis in $L^2_{\text{per}}(\br^2/\Lambda)$. The eigenvalues $E_b(\bk)$ are called band dispersion functions or Bloch bands which are Lipschitz continuous, and $E_b(\bk)=0$ if and only if $b=1$ and $\bk=\boldsymbol{0}$ with the corresponding normalized eigenfunction $\Phi_1(\bx;\boldsymbol{0})=\phi_1(\bx;\boldsymbol{0})\equiv |\Omega|^{-\frac12}$. Thus, for any $\lambda>0$, there exist a positive constant $C$ such that
\begin{equation} \label{eigenvalue:lowbd}
	E_b(\bk)\ge C>0,\quad\text{for}~|\bk|\ge \dlt_{b,1}\lambda.
\end{equation}
We will choose an appropriate constant $\lambda$ for the proof convenience in Section \ref{proof:prop}.
%After giving the Dirac points $(\gk_*,~ E_D)$ in the next subsection, we can choose $B=|\gk_*|$.

The corresponding quasi-periodic eigenfunctions $\Phi_b(\bx;\bk)$ are called Bloch modes. For any given $\bk\in\Omega^*$, $\Phi_{b}(\bx;\bk)~(\text{or}~\phi_b(\bx;\bk)),~b\ge1$ is analytic for $\bx \in \br^2$ from the regularity theory of elliptic operator. Moreover, the set of all Bloch modes  $\{\Phi_{b}(\bx; \bk)\}_{b\ge1,~
\bk\in\Omega^*}$ is complete in $L^2(\br^2)$. That is for any $f\in L^2(\br^2)$,
\begin{eqnarray}
	f(\bx) = \frac{1}{|\Omega^*|}\sum_{b\ge1}\int_{\Omega^*}\tilde{f}_b(\bk)\Phi_b(\bx;\bk)~ d\bk, \quad \text{where} \quad \tilde{f}_b(\bk) = \bla\Phi_{b}(\cdot;\bk),~f(\cdot)\bra, \label{fb:f}
\end{eqnarray}
and the Parseval formula for the Bloch decomposition in $L^2(\br^2)$ holds,
\begin{equation}
	\big\|f\big\|^2_{L^2(\br^2)}=\frac{1}{|\Omega^*|}\sum_{b\geq1}\int_{\Omega^*}|\tilde{f}_b(\bk)|^2~d\bk~. \label{l2norm}
\end{equation}
Thanks to the Weyl's law, i.e., $E_b(\bk)\approx b~(b\gg 1)$ uniformly for all $\bk\in\Omega^*$, and then any $f\in H^s(\br^2)$, $\|f\|_{H^s(\br^2)}$ can be approximated by
\begin{align}
\nonumber	\|f\|^2_{H^s(\br^2)} \approx&~ \bla \big(1+\mathcal{L}^A\big)^s f(\bx), f(\bx)\bra \\
\nonumber	=&~\frac{1}{|\Omega^*|} \sum_{b\ge 1}\int_{\Omega^*}(1+E_b(\bk))^s|\tilde{f}_b(\bk)|^2~d\bk \\
\approx&~ \sum_{b\ge 1}(1+b)^s\int_{\Omega^*}|\tilde{f}_b(\bk)|^2~d\bk . \label{sobolevnorm}
\end{align}

\subsection{Honeycomb structured material weight and Dirac points}

Before introducing the honeycomb structured material weight, we define the following symmetry operators acting on a function $f(\bx)$ defined in $\br^2$.
Parity inversion operator $\mathcal{P}$: $\mathcal{P}[f](\bx)=f(-\bx)$;
Complex conjugate operator $\mathcal{C}$: $\mathcal{C}[f](\bx)=\overline{f(\bx)}$;
$\frac{2\pi}{3}-$rotation operator $\mathcal{R}$: $\mathcal{R}[f](\bx)=f(R^*\bx)$
where $R$ is $\frac{2\pi}{3}-$clockwise rotation matrix
\begin{equation}
	R=
	\begin{pmatrix}
	-\frac12 & \frac{\sqrt{3}}{2} \\
	-\frac{\sqrt{3}}{2} & -\frac12
	\end{pmatrix}.
\end{equation}

A honeycomb structured material weight is defined as follows
\begin{definition}\label{hcdef}
A $2\times2$ complex-valued matrix function $A(\bx)$ is called a honeycomb structured material weight if it satisfies the Definition \ref{assume} and further
	\begin{enumerate}
	\setcounter{enumi}{0}
	\item $A(\bx)$ is $\Lambda$-periodic, i.e., $A(\bx+\bv)=A(\bx),~\forall~\bx\in\br^2,~\bk\in\Lambda$;
	\item $A(\bx)$ is $\mathcal{PC}$-invariant, i.e., $\mathcal{PC}[A](\bx)\equiv\overline{A(-\bx)}=A(\bx),~\forall~\bx\in\br^2$;
		\item $A(\bx)$ satisfies $\mathcal{R}[A](\bx)\equiv A(R^*\bx)=R^*A(\bx)R,~\forall~\bx\in\br^2$.
	\end{enumerate}
\end{definition}
A consequence of $A(\bx)$ being a honeycomb structured material weight is that $\mathcal{PC}\ml^A=\ml^A\mathcal{PC},~ \mathcal{R}\ml^A=\ml^A\mathcal{R}$. Generically, this leads to the existence of Dirac points in the dispersion surfaces of $\ml^A=-\nabla\cdot \big(A(\bx)\nabla\big)$. The specific definition of Dirac points is given as follows, see \cite{fefferman2012honeycomb,lee2018elliptic}.

\begin{definition} \label{Dirac pts}
The quasi-momentum/eigenvalue pair $(\gk_*, E_D)\in\Omega^*\times \br_+$ is called a Dirac point if there exists an integer $b_*\ge 1$ and Floquet-Bloch eigenpairs mappings
\begin{equation*}
    \bk\mapsto(\Phi_{b_*}(\bx;\bk),~E_{b_*}(\bk))~~\text{and}~~\bk\mapsto(\Phi_{b_*+1}(\bx;\bk),~E_{b_*+1}(\bk))
\end{equation*}
such that:
\begin{enumerate}
\item[(1)]  $E_D=E_{b_*}(\gk_*)=E_{b_*+1}(\gk_*)$ is a two-fold degenerate $L^2_{\gk_*}$-- eigenvalue of $\mathcal{L}^A$. There exists two orthogonal eigenfunctions $\Phi_1(\bx),~\Phi_2(\bx)=\overline{\Phi_1(-\bx)}$,
\begin{eqnarray}
\bla\Phi_i(\bx),~\Phi_j(\bx)\bra_{L^2(\Omega)}=\delta_{ij},\quad i,~j=1,~2; \label{Phi:ij}
\end{eqnarray}
\item[(2)] Denote that $E_{-}(\bk)=E_{b_*}(\bk),~E_{+}(\bk)=E_{b_*+1}(\bk)$, and $\bkappa=(\kappa_1,\kappa_2)=\bk-\gk_*$.  There exist  $v_{_F}>0$ and $q_0>0$, such that for $0\le |\bkappa| \le q_0$,
\begin{equation}
   E_{\pm}(\gk_*+\bkappa)-E_D = \pm v_{_F}|\bkappa|(1+e_{\pm}(\bkappa)),
\end{equation}
where  $|e_{\pm}(\bkappa)|\le C|\bkappa|$ for some constant $C$.
\end{enumerate}
\end{definition}

Consider the two high symmetric points in $\Omega^*$, $\gk=\frac{1}{3}(\bk_1-\bk_2)$ and $\gk'=-\gk$. In the previous literature, Lee-Thorp, Weinstein and Zhu proved, if $A(\bx)$ is a honeycomb structured material weight in the sense of Definition \ref{hcdef}, the Dirac points generically appear in two adjacent dispersion surfaces of $\ml^A$ conically intersecting at $\gk_*=\gk$ and $\gk'$, see Theorem 4.2 and 4.10 in \cite{lee2018elliptic}.

Furthermore, we have two more conclusions at the following proposition
\begin{proposition}\label{dirac:modes}
	Let $(\gk_*, E_D)$ be a Dirac point in the sense of Definition \ref{Dirac pts}. If $b\notin \{+,-\}$, there exists a constant $C>0$, and $q_1>0$ small, such that for any $0\le |\bk-\gk_*| \le q_1$
	\begin{equation}\label{lowerbound}
	|E_{b}(\bk)-E_D|\ge C.
	\end{equation}
	Let $\Phi_{-}(\bx;\bk)=\Phi_{b_*}(\bx;\bk), ~\Phi_+(\bx;\bk)=\Phi_{b_*+1}(\bx;\bk)$. When
	$0<|\bkappa|\le q_0$, one can expand $\Phi_{\pm}(\bx;\bk)$ in the following form,
	\begin{equation}\label{mode:expan}
	\Phi_{\pm}(\bx;\gk_*+\bkappa) = \frac{e^{i\bkappa\cdot\bx}}{\sqrt{2}}\Big[ \frac{\kappa_1+i\kappa_2}{|\bkappa|}\Phi_1(\bx)\pm \Phi_2(\bx)+\mathcal{O}_{H^2_{\gk_*}(\br^2/\Lambda)}(|\bkappa|)\Big].
	\end{equation}
\end{proposition}

The proof of \eqref{lowerbound} in Proposition \ref{dirac:modes} is a direct consequence of \eqref{sequence:eigenva} and the Lipschitz continuity for all eigenvalues, while for the rigorous proof of \eqref{mode:expan} is referred to Theorem 3.2 in \cite{fefferman2014wave}.

In this work, we consider a slowly modulated honeycomb structured material weight $W_{\veps}(\bx)=A(\bx)+\veps \kappa(\veps \bx)B(\bx)$. Throughout this work, we require the following assumptions on $W_{\veps}(\bx)$ hold

\begin{assum}\label{assumption1}
	The $2\times2$ matrix $W_{\veps}(\bx)=A(\bx)+\veps \kappa(\veps \bx)B(\bx)$ is a material weight in the sense of Definition \ref{assume}, and
	\begin{enumerate}[label=(\subscript{A}{\arabic*}), itemindent=1em]
	\item $A(\bx)$ is a honeycomb structured material weight in the sense of Definition \ref{hcdef};
	\item $\kappa(\gx)$ is a real scalar function in $ C^{\infty}_b(\br^2)$;
	\item $B(\bx)$ is smooth, Hermitian, $\Lambda$-periodic, and $B(-\bx)=B(\bx)$, $\overline{B(\bx)}=-B(\bx)$;
    %\item $\delta$ is a sufficiently small positive parameter.
\end{enumerate}
\end{assum}

\noindent\textbf{Remark:} Assumption $(A_3)$ implies that $\mathcal{PC}\mathcal{L}^B=-\ml^B\mathcal{PC}$ and further $B(\bx)=b(\bx)\sigma_2$ with $b(\bx)$ real and even. This assumption is related to the time-reversal symmetry breaking in the study of photonic topological insulators in real applications \cite{haldane2008Possible,khanikaev2013photonic,lee2018elliptic}. It also indicates that the second order operator $\mlb$ is actually a first order operator on $\eta$, i.e.,
\begin{equation}\label{firstorder}
	-\nabla\cdot \big(\kappa(\veps \bx)B(\bx)\nabla\eta\big) = -\nabla\cdot\big[\nabla\cdot\big(\kappa(\veps \bx)B(\bx)\big)\eta\big].
\end{equation}
In practical applications, there is another way to break $\mathcal{PC}$-symmetry by assuming $B(\bx)=\tilde{b}(\bx)I_{2\times2}$ with $\tilde{b}(\bx)$ real and odd, which is related to parity symmetry breaking \cite{lee2018elliptic}. In this case $\mlb$ remains a second order operator and our analysis does not apply.

A typical example of a media satisfying Assumption \ref{assumption1} in real applications is the magneto-optical material given in Haldane and Ragu's work \cite{haldane2008Possible}. In their physical setup, the material weights written in our notation are
\begin{equation*}
	A(\bx)=\epsilon^{-1}(\bx) \textrm{I}_{2\times2},~~~B(\bx)=\epsilon^{-2}(\bx)\sigma_2,~~~\veps\kappa(\veps \bx)=\gamma(\veps \bx),
\end{equation*}
where $\epsilon(\bx)$ is the electric permittivity which is even, real and $\mathcal{R}$-invariant, and $\gamma(\veps\bx)$ represents the Fardy-rotation effect which is assumed to be small and slowly varying.
%and we have
%\begin{equation*}
%\ptl_t^2 H_3-\nabla_{\perp}\cdot\big(\epsilon^{-1} \textrm{I} ~\nabla_{\perp}H_3\big)-\nabla_{\perp}\cdot\big(\gamma\epsilon^{-2}\sigma_2\nabla_{\perp}H_3\big)=0
%\end{equation*}

We end this section by including the following important results \cite{lee2018elliptic}.
\begin{proposition}\label{propinner}
	Suppose that $\Phi_1(\bx),~\Phi_2(\bx)=\overline{\Phi_1(-\bx)}$ are eigenfunctions of $\ml^A$ with respect to the Dirac point $(E_D,\gk_*)$ given in Definition \ref{Dirac pts}, $B(\bx)$ satisfies $(A_3)$, and define
\[
	\mathcal{A} = \frac1i \Big(A(\bx)\nabla+\nabla\cdot A(\bx) \Big),\quad \mathcal{A}\Phi_j=\frac1i \Big(A(\bx)\nabla\Phi_j+\nabla\cdot \big(A(\bx)\Phi_j\big) \Big),\quad j=1,~2.
\]
Then, the following identities hold:
	\begin{equation}\label{iA:12}
	\begin{aligned}
		&\bla \Phi_1(\bx),~\mathcal{A}\Phi_2(\bx)\bra_{L^2(\Omega)}=v_{_F}
		\begin{pmatrix}
			1 \\
			i
		\end{pmatrix},\quad
		\bla \Phi_2(\bx),~\mathcal{A}\Phi_1(\bx)\bra_{L^2(\Omega)}=v_{_F}
		\begin{pmatrix}
			1 \\
			-i
		\end{pmatrix}, \\
		&\bla \Phi_1(\bx),~ \mathcal{A}\Phi_1(\bx)\bra_{L^2(\Omega)}=\bla \Phi_2(\bx),~ \mathcal{A}\Phi_2(\bx)\bra_{L^2(\Omega)}\equiv 0,
	\end{aligned}
	\end{equation}
	and
	\begin{equation}\label{lb:11}
	\begin{aligned}
		&\bla \Phi_1(\bx),~\ml^B\Phi_1(\bx) \bra_{L^2(\Omega)} =-\bla \Phi_2(\bx),~\ml^B\Phi_2(\bx) \bra_{L^2(\Omega)} , \\
		&\bla \Phi_1(\bx),~\ml^B\Phi_2(\bx)\bra_{L^2(\Omega)} = \bla \Phi_2(\bx),~\ml^B\Phi_1(\bx) \bra_{L^2(\Omega)} = 0.
	\end{aligned}
	\end{equation}
\end{proposition}
Here we define $\vartheta_{\sharp}=\bla \Phi_1,~\ml^B\Phi_1 \bra_{L^2(\Omega)}$ and assume $\vartheta_{\sharp}\neq 0$ in this work.

\section{Dirac equation with a varying mass}
This paper is to show that the Dirac equation \eqref{dirac} governs the envelope dynamics of the wave problem \eqref{weq:perturbed} under a prescribed initial condition \eqref{ini}. In this section, we first state results on well-posedness, estimates on solutions to the Dirac equation \eqref{dirac} in the following proposition.
\begin{proposition}\label{schwartz}
Let $E_D>0$, $v_{_F}>0$, $\vartheta_{\sharp}\neq 0$ be given constants as before, $\kappa(\gx)\in C_b^{\infty}(\br^2)$, and $\balpha_0(\gx)=\big(\alpha_{10}(\gx),\alpha_{20}(\gx)\big)\in\mathcal{S}(\br^2)$. Then, for any $0<\rho<+\infty,~s>2$, the Dirac equation \eqref{dirac} has a unique solution $\balpha(\gx, T)=\big(\alpha_1(\gx, T),~\alpha_2(\gx, T)\big)^T \in C^{\infty}\big([0,\infty)\times\br^2\big)$  and
\begin{equation}\label{dirac:hsnorm}
\balpha(\gx,T)\in C^0\big([0,\rho], H^{s}(\br^2)\big)\cap C^1\big([0,\rho], H^{s-1}(\br^2)\big).
\end{equation}\label{wellpose}
Moreover, for any $T\ge0,~l\in \mathbb{N}$, $\ptl_T^l\alpha_j(\cdot, T)\in \mathcal{S}(\br^2)$, specifically, $\forall ~M\in \mathbb{N}$, $\mathbf{n}\in \mathbb{N}^2$, there exists a constant $C>0$ such that
\begin{equation} \label{schwartzestimate}
	\sup_{T\in[0,\rho],~\gx\in \br^2}\bav(1+|\gx|^2)^{\frac M2}\ptl_{\gx}^{\mathbf{n}}\ptl_T^l \balpha(\gx, T)\bav < C.
\end{equation}
\end{proposition}

Note that Dirac equation \eqref{dirac} is actually a first order linear hyperbolic system. If $\kappa(\gx)$ is a constant, \eqref{dirac} has a unique solution in Schwartz space for the time $T\in[0,+\infty)$ in term of the Fourier transform arguments \cite{fefferman2014wave}. However, for a general $\kappa(\gx)\in C_b^{\infty}(\br^2)$, a comprehensive proof to the above proposition will be postponed in the Appendix B.
\newline

One of the most interesting applications of the reduced Dirac equation \eqref{dirac} is its capability to describe dynamics of topological edge states. To be more specific, suppose that $\kappa(\gx)$ is a domain wall function, i.e., $\kappa(\gx)=\tilde{\kappa}(\zeta)$ with $\tilde\kappa(\zeta)\to \pm \kappa_\infty$ as $\zeta\to\pm\infty$, where $\zeta=\boldsymbol{\mathfrak{K}}\cdot\gx\in\br$ and $\boldsymbol{\mathfrak{K}}=(\mathfrak{K}_1,\mathfrak{K}_2)\in \Lambda^*$, $\kappa_\infty>0$, see \cite{fefferman2017topologically,lee2018elliptic} for details. Hereafter, we drop the tilde on top of $\kappa$.

Let
\begin{align*}
\boldsymbol{\mathfrak{K}}^{\perp}=(-\mathfrak{K}_2,\mathfrak{K}_1),	\quad \xi=\boldsymbol{\frk}^{\perp}\cdot\gx.
\end{align*}

Rewriting the Dirac equation \eqref{dirac} in new coordinates system $(\xi,\zeta)$ yields that
\begin{equation}\label{newcoor}
\left\{
	\begin{aligned}
		i\ptl_{T}\alpha_1+\frac{v_{_F}}{2\sqrt{E_D}}\big(\frk_1+i\frk_2\big)\ptl_{\xi}\alpha_2-\frac{v_{_F}}{2\sqrt{E_D}}\big(i\frk_1-\frk_2\big)\ptl_{\zeta}\alpha_2+\frac{\vartheta_{\sharp}\kappa(\zeta)}{2\sqrt{E_D}}\alpha_1=0, \\
		i\ptl_{T}\alpha_2-\frac{v_{_F}}{2\sqrt{E_D}}\big(\frk_1-i\frk_2\big)\ptl_{\xi}\alpha_1-\frac{v_{_F}}{2\sqrt{E_D}}\big(i\frk_1+\frk_2\big)\ptl_{\zeta}\alpha_1-\frac{\vartheta_{\sharp}\kappa(\zeta)}{2\sqrt{E_D}}\alpha_2=0,
	\end{aligned}
\right.
\end{equation}
where $\alpha_j=\tilde{\alpha}_j(\xi,\zeta, T)=\alpha_j(\bx(\xi,\zeta), T)$ and we have used the same symbols before and after changing coordinates for notational convenience.

We are interested in a particular solution to \eqref{newcoor} which decays to zero as $|\zeta|\to \infty$ and keeps periodic in $\xi-$direction. This solution is referred to as the topological edge state. Namely, let
\begin{equation*}
	\balpha(\xi, \zeta, T)=e^{ik_\parallel\xi-i\mu(k_\parallel)T}\bbeta(\zeta; k_\parallel).
\end{equation*}
It deduces that $\bbeta(\zeta; k_\parallel)$ satisfies the following eigenvalue problem in $L^2(\mathbb{R})$,
\begin{equation}\label{edge}
	\mathcal{D}(k_\parallel)\bbeta(\zeta;k_\parallel)=\mu(k_\parallel)\bbeta(\zeta;k_\parallel),
\end{equation}
where $\mathcal{D}(k_\parallel)$ is the 1D Dirac operator
\begin{equation*}
\mathcal{D}(k_\parallel)= \frac{v_{_F}}{2\sqrt{E_D}}
\begin{pmatrix}
0 & i\mathfrak{K}_1-\mathfrak{K}_2 \\
i\mathfrak{K}_1+\mathfrak{K}_2 & 0
\end{pmatrix}\ptl_{\zeta}
-\frac{v_{_F}}{2\sqrt{E_D}}k_\parallel
\begin{pmatrix}
0 & i\mathfrak{K}_1-\mathfrak{K}_2 \\
-i\mathfrak{K}_1-\mathfrak{K}_2 & 0
\end{pmatrix}
-\frac{\vartheta_\sharp\kappa(\zeta)}{2\sqrt{E_D}}\sigma_3.
\end{equation*}

In the previous work, Lee-Thorp, Weinstein and Zhu \cite{lee2018elliptic} derived the same equation as \eqref{edge} in the case that $k_\parallel=0$. In this scenario, there exist the so-called zero-energy states for $\mathcal{D}(0)$. That is, $\big(\bbeta(\zeta;0),0\big)$ is the eigenpair of the operator $\mathcal{D}(0)$ in $L^2(\mathbb{R})$ with
\begin{equation}
	\bbeta(\zeta;0)=\left\{
	\begin{aligned}
		\frac{\sqrt{2}}2\gamma e^{-\frac{|\vartheta_{\sharp}|}{v_{_F}|\boldsymbol{\mathfrak{K}}|}\int_0^{\zeta}\kappa(s)ds}
		\begin{pmatrix}
		\frac{-\mathfrak{K}_2+i\mathfrak{K}_1}{|\boldsymbol{\mathfrak{K}}|} \\
		-1
		\end{pmatrix},\quad \text{if}\quad\vartheta_{\sharp}>0,\\
		\frac{\sqrt{2}}2\gamma e^{-\frac{|\vartheta_{\sharp}|}{v_{_F}|\boldsymbol{\mathfrak{K}}|}\int_0^{\zeta}\kappa(s)ds}		
		\begin{pmatrix}
		\frac{-\mathfrak{K}_2+i\mathfrak{K}_1}{|\boldsymbol{\mathfrak{K}}|} \\
		1
		\end{pmatrix},\quad \text{if}\quad\vartheta_{\sharp}<0,
	\end{aligned}
	\right.
\end{equation}
and here $\gamma$ is the normalization constant.

Our derivation and analysis demonstrate the existence of edge states proved in \cite{lee2018elliptic} from the evolutionary effective envelope equation \eqref{dirac}. Equation \eqref{dirac} can be further used to study the dynamics of such states as well as their interactions with defects, perturbations by manipulating $\kappa(\zeta, \xi)$.
\begin{figure}[!hbt]
\centering
\begin{subfigure}{0.2\textwidth}
	\includegraphics[width=\textwidth]{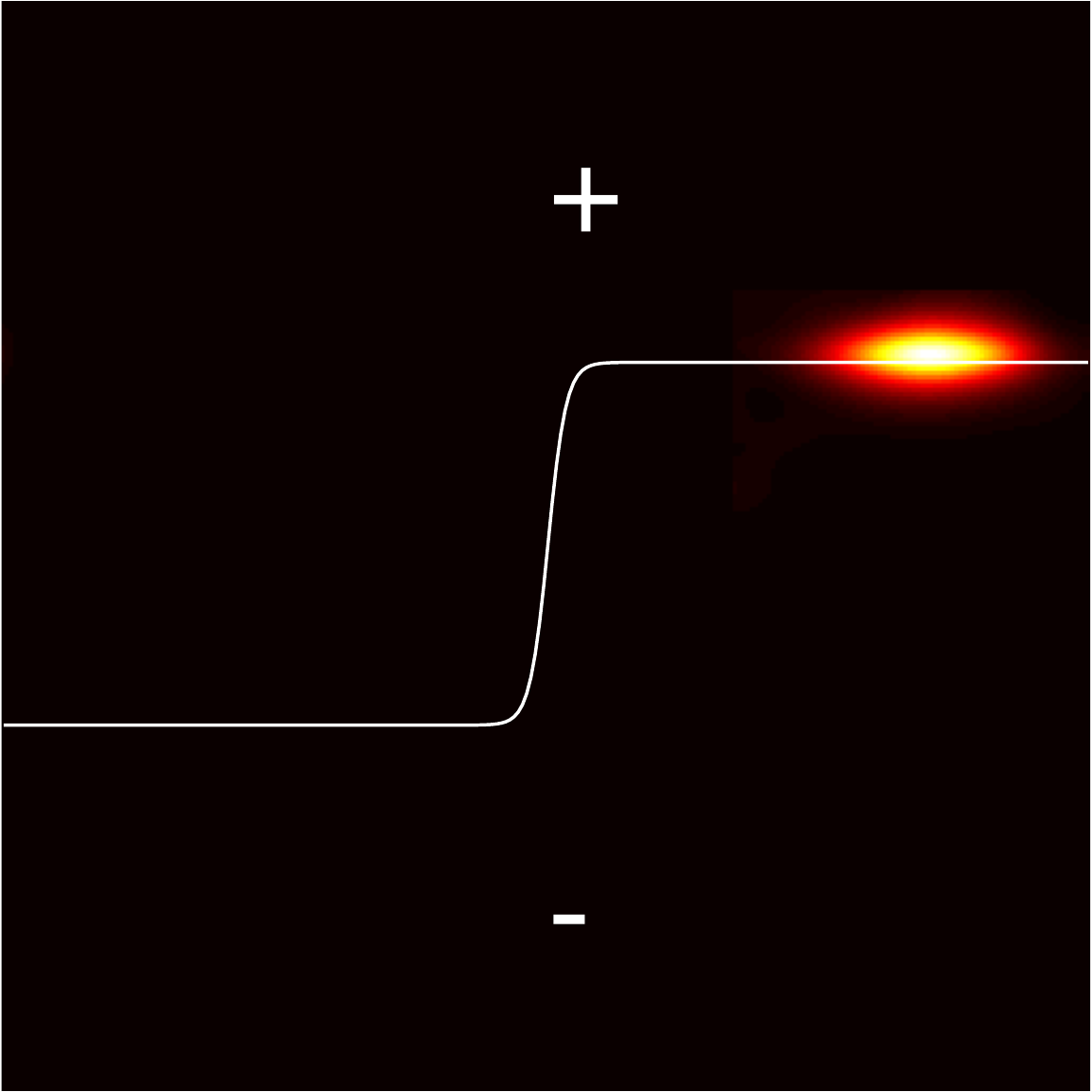}
\end{subfigure}\quad
\begin{subfigure}{0.2\textwidth}
	\includegraphics[width=\textwidth]{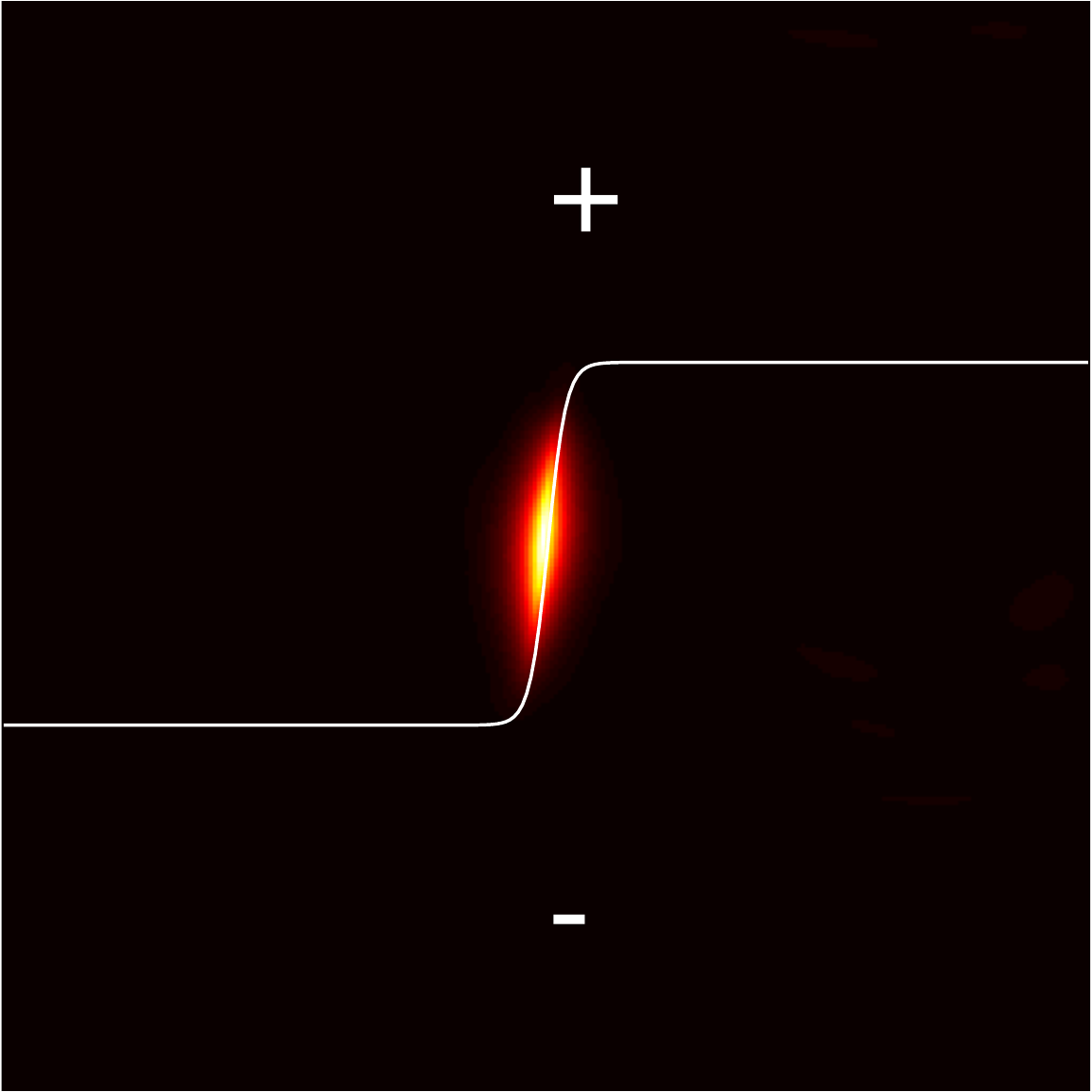}
\end{subfigure}\quad
\begin{subfigure}[a]{0.2\textwidth}
	\includegraphics[width=\textwidth]{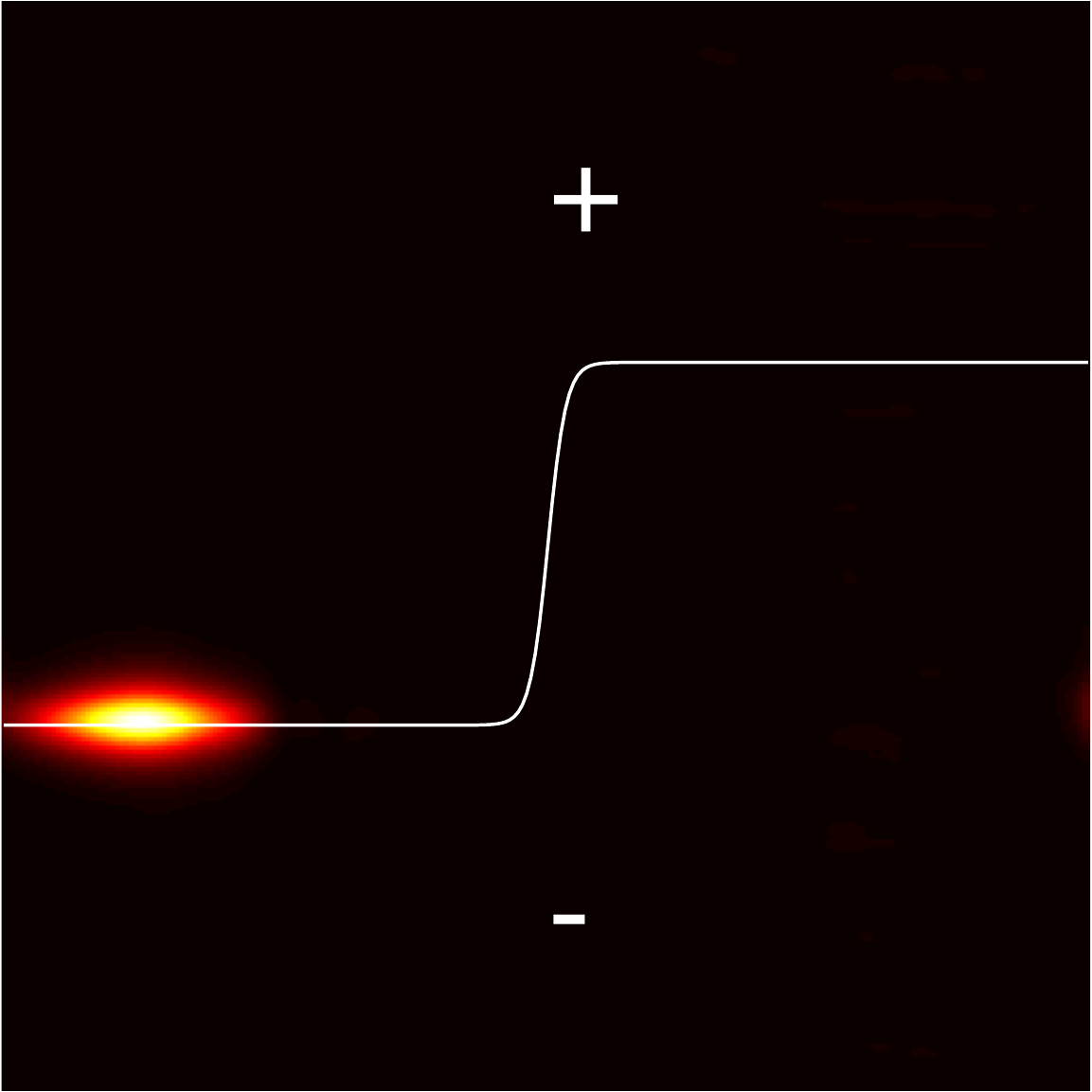}
\end{subfigure}

\begin{subfigure}{0.2\textwidth}
	\includegraphics[width=\textwidth]{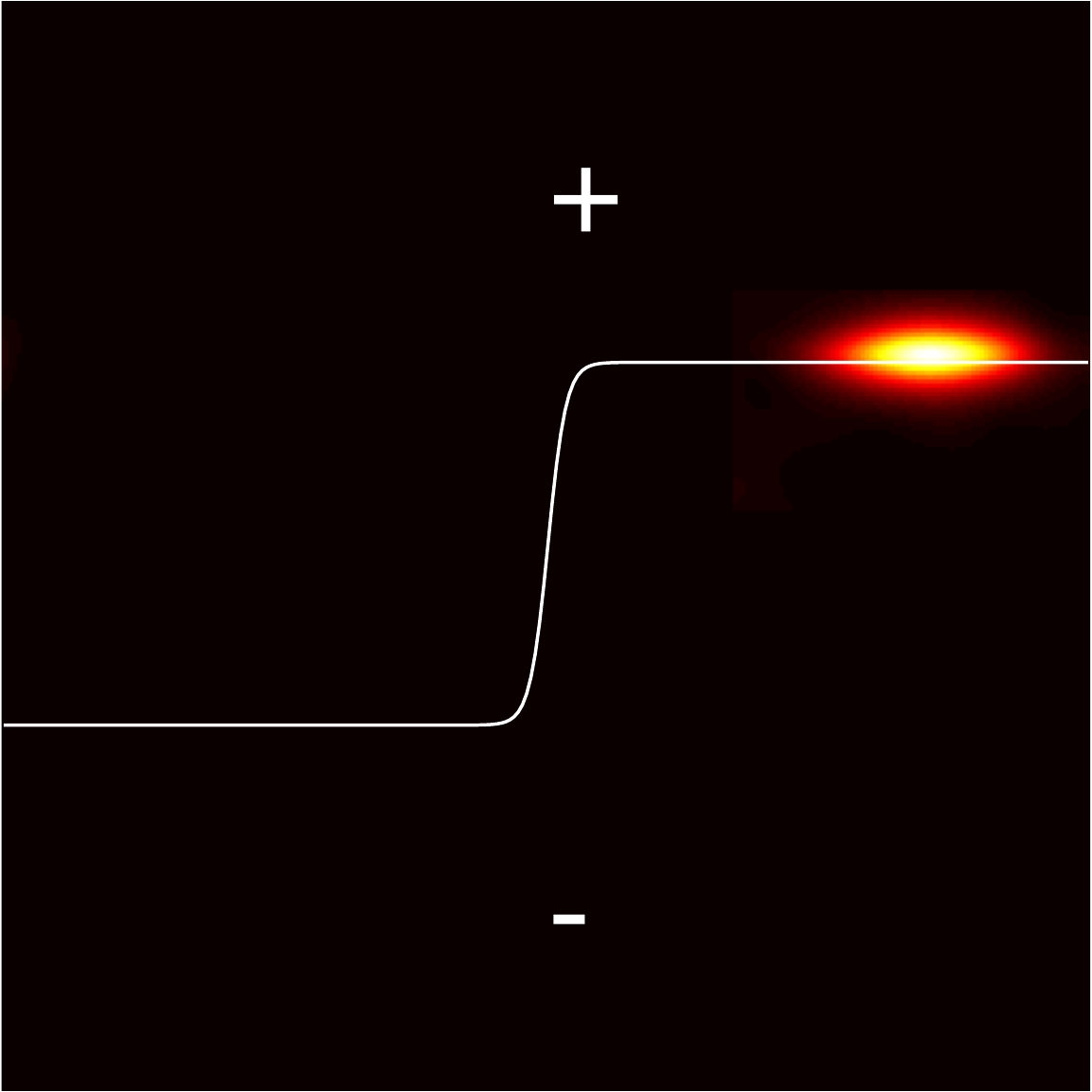}
	\caption{T=0}
\end{subfigure}\quad
\begin{subfigure}{0.2\textwidth}
\includegraphics[width=\textwidth]{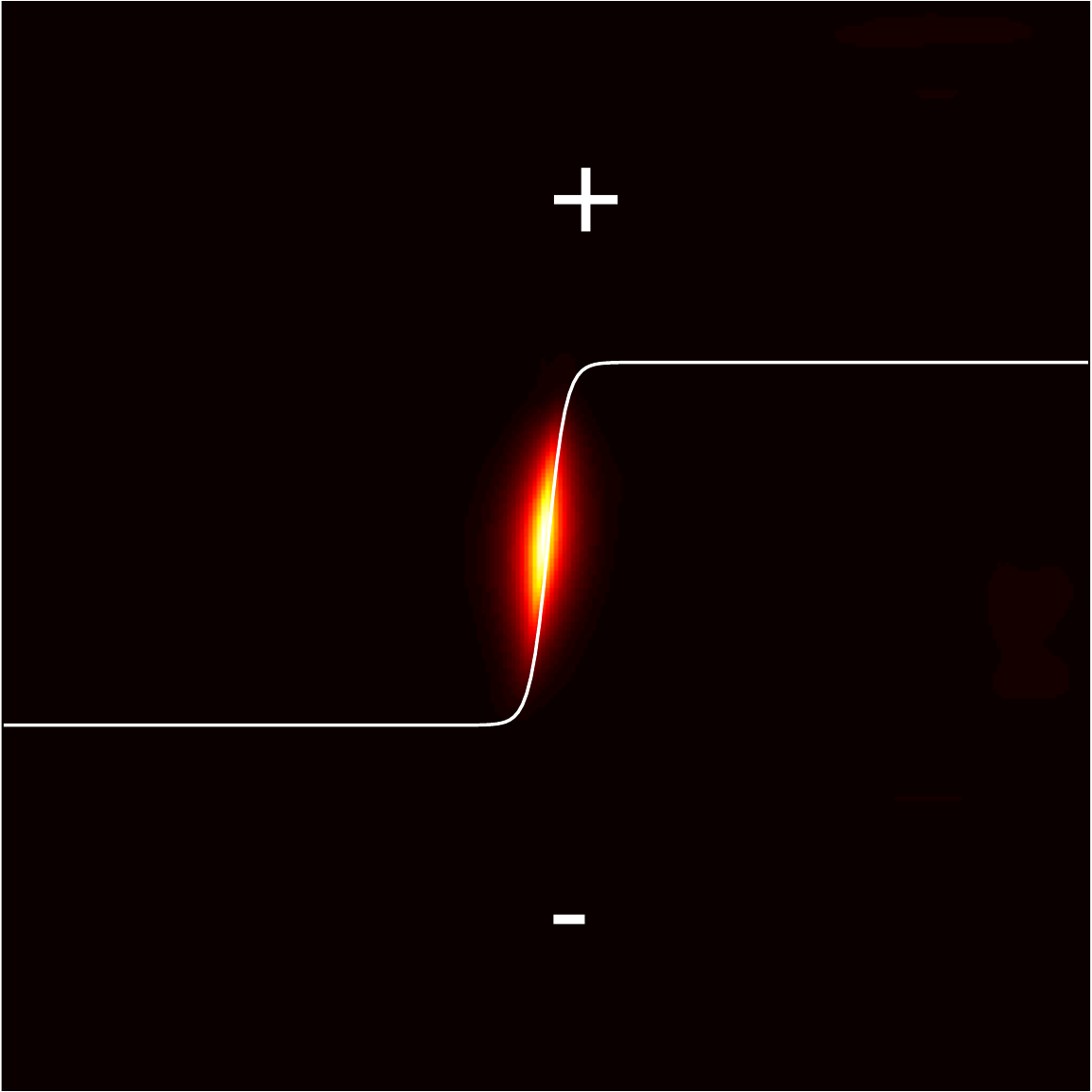}
\caption{T=30}
\end{subfigure}\quad
\begin{subfigure}{0.2\textwidth}
\includegraphics[width=\textwidth]{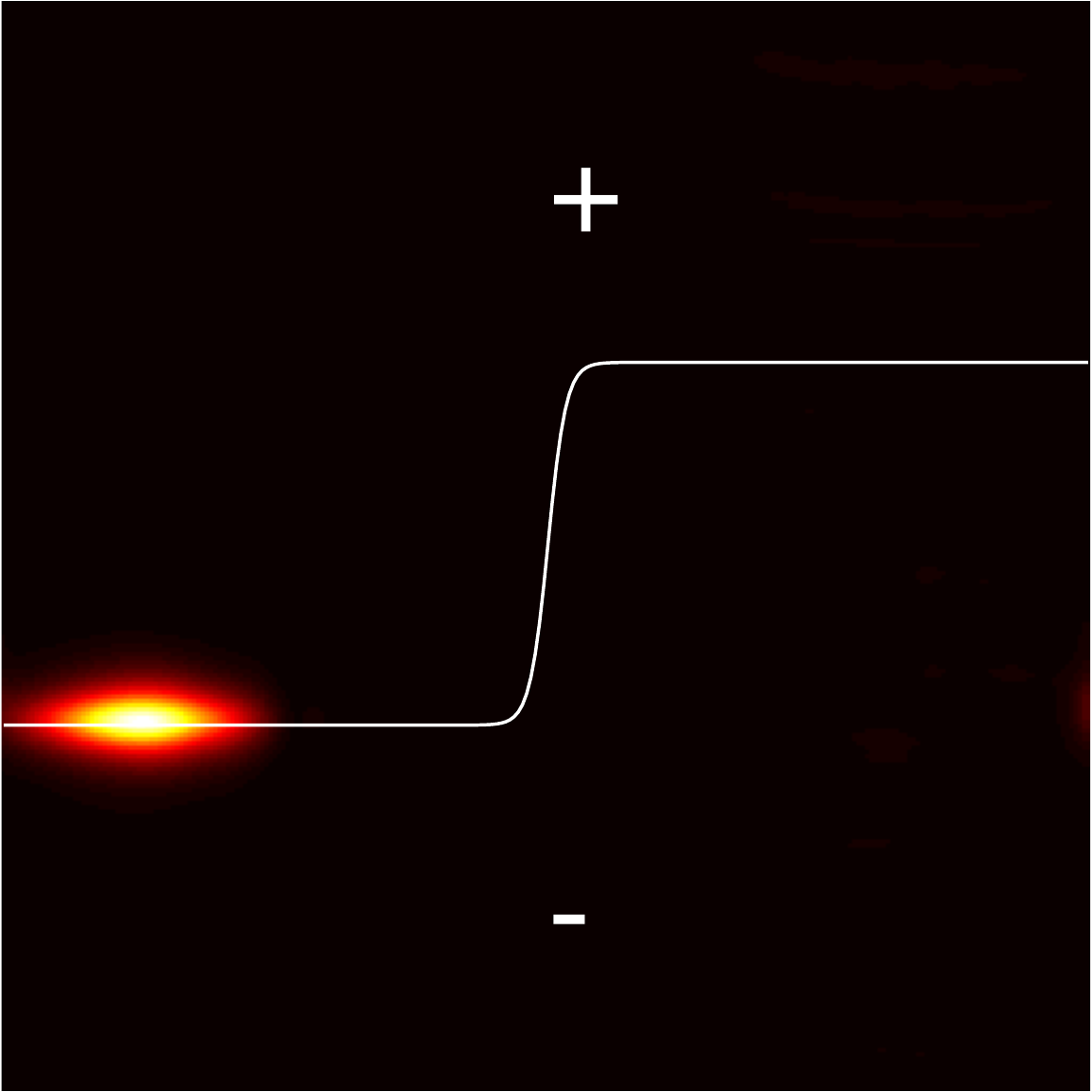}
\caption{T=60}
\end{subfigure}
\caption{The plots of the intensity of the solution to equation \eqref{dirac} with initial condition \eqref{dirac:ini} at different times. Top panel: $|\alpha_1|$. Bottom panel: $|\alpha_2|$. The white curve represents the ``edge", i.e.,  above the curve $\kappa(\gx)$ is positive, below the curve $\kappa(\gx)$ is negative, and $\kappa(\gx)$ vanishes on the curve. }
\end{figure}
	
In order to show that the envelope equation \eqref{dirac} can exhibit many interesting solutions which describe the novel and subtle physical phenomena, we illustrate a typical propagation pattern by solving \eqref{dirac} numerically. In this simulation, the coefficients of the equation is normalized for simplicity. Alternatively, we set $\frac{v_{_F}}{2\sqrt{E_D}}=1,~ \frac{\vartheta_{\sharp}}{2\sqrt{E_D}}=1$. The modulation that we choose is
\begin{equation}
\kappa(\gx)=\tanh(X_2-10\tanh(X_1)).
\end{equation}
The initial condition is
\begin{equation}
\alpha_{10}(\gx)=\text{sech}(X_2-\textsf{X}_2^0)e^{-(X_1-\textsf{X}_1^0)^2}~\text{with}~ \textsf{X}_2^0=10\tanh(\textsf{X}_1^0),~\text{and}~ \alpha_{20}=-\alpha_{10}.
\end{equation}

In Figure 1, we plot several snapshots of the intensity of the solution to equation \eqref{dirac} at three successive times. It can be seen that the waves travel along the edge without any energy leaking to the bulk or traveling back. This interesting phenomenon is related to the topologically protected wave propagation which is one of the current focuses in many applied fields. Our rigorous justification of \eqref{dirac} from \eqref{weq:perturbed} provides a solid mathematical foundation for such interesting problems. Due to the length and scope of this paper, we leave the further analysis on the reduced envelope equation \eqref{dirac} in future works.

\section{Main results}

The main goal of this paper is to show that the equation \eqref{weq:perturbed} with initial condition \eqref{ini} has the asymptotic solution \eqref{solu1} with the envelopes $\alpha_j(\veps \bx,\veps t),~ j=1,~2$ satisfying the Dirac equation \eqref{dirac}\eqref{dirac:ini}. The well-posedness of  Cauchy problem \eqref{weq:perturbed}\eqref{ini} is a standard result by the theory on linear hyperbolic systems, see e.g., \cite{kato1975cauchy,racke1992lectures}. Our task is reduced to rigorously justify that the error $\eta(\bx,t)$ is small over a large but finite time. To this end, we substitute \eqref{solu1} into \eqref{weq:perturbed}\eqref{ini} and obtain the equation of $\eta(\bx, t)$,
\begin{eqnarray}
   \partial^2_{t}\eta+\mathcal{L}^A\eta+\veps\mlb\eta = e^{i\sqrt{E_D}t}\big(F_1(\bx,t)+F_2(\bx,t)\big),
\end{eqnarray}
with the initial condition
\begin{equation}
\begin{split}
\eta(\bx,0) &= 0,\\
\partial_t\eta(\bx,t)|_{t=0} &= -\veps^2 \partial_T\alpha_j(\veps \bx,0)\Phi_j(\bx)	:= F_0(\bx),
\end{split}
\end{equation}
where $F_1(\bx,t)$ and $F_2(\bx, t)$ are
\begin{align}\label{F1}
\nonumber  F_1(\bx,t) =&~ -\veps^2 \Big[2i\sqrt{E_D}\ptl_T\alpha_j\Phi_j(\bx)-\nabla_{\gx}\alpha_j\cdot\big[A(\bx)\nabla\Phi_j(\bx)+\nabla\cdot\big(A(\bx)\Phi_j(\bx)\big)\big] \\
\nonumber  &~ -\kappa(\veps\bx)\alpha_j\nabla\cdot\big(B(\bx)\nabla\Phi_j(\bx)\big)\Big] \\
	=&~ -\veps^2 \Big[2i\sqrt{E_D}\ptl_T\alpha_j\Phi_j(\bx)-\nabla_{\gx}\alpha_j\cdot i\mathcal{A}\Phi_j(\bx)+\kappa(\veps\bx)\alpha_j\ml^B\Phi_j(\bx)\Big],
\end{align}
\begin{align}\label{F2}
\nonumber  F_2(\bx,t) =&~ -\veps^3 \Big[\ptl_T^2\alpha_j\Phi_j(\bx) -\kappa(\veps\bx)\nabla_{\gx}\alpha_j\cdot\big[B(\bx)\nabla\Phi_j(\bx)+\nabla\cdot\big(B(\bx)\Phi_j(\bx)\big)\big] \\
\nonumber  &~ -\big(\nabla_{\gx}\kappa(\veps \bx)\alpha_j\big)\cdot \big(B(\bx)\nabla\Phi_j(\bx)\big)-\nabla_{\gx}^2\alpha_j:A(\bx)\Phi_j(\bx) \Big] \\
  &~ +\veps^4 \Big[\kappa(\veps \bx)\nabla_{\gx}^2\alpha_j:B(\bx)\Phi_j(\bx)-\big(\nabla_{\gx}\kappa(\veps \bx)\nabla_{\gx}\alpha_j\big):B(\bx)\Phi_j(\bx)\Big].
\end{align}

The main result is concluded as follows.
\begin{theorem}\label{mainthm}
Suppose $W_{\veps}(\bx)$ satisfies Assumption \ref{assumption1}, the operator $\ml^A=-\nabla\cdot \big(A(\bx)\nabla\big)$ has a Dirac point $(\gk_*, E_D)$,  $\Phi_1(\bx)$ and $\Phi_2(\bx)$ are the associated eigenfunctions given in Definition \ref{Dirac pts},  and $\alpha_{10}(\gx)$, $\alpha_{20}(\gx)$ are Schwartz functions of $\gx\in\br^2$. Then the wave equation \eqref{weq:perturbed} with the initial condition \eqref{ini} has a unique solution of the form \eqref{solu1}, where $\alpha_j(\gx,T),~j=1,~2$ are the solution of the system \eqref{dirac}\eqref{dirac:ini}, and for any $s\ge 0$, $\rho>0$, $0<\nu<1$ and $\veps>0$ sufficiently small,
\begin{equation}\label{mainthmresult}
\sup_{0\le t\le\rho\veps^{^{-1}}}\|\eta(\bx, t)\|_{H^s(\br^2)} \le C\veps^{1-\nu}.
\end{equation}
Here $C$ is independent of $\veps$.
\end{theorem}

\noindent{\bf Remark:} In the literatures on wave packet problems, the macroscopic reference frame $\gx=\veps\bx,~T=\veps t$ with a $\veps-$scaled Sobolev space $H^s_{\veps}(\br^2)$, i.e., $\forall~f(\gx)\in H^s_{\veps}(\br^2)$,
\begin{equation*}
	\|f(\gx)\|^2_{H^s_{\veps}(\br^2)}=\sum_{|\gamma|\le s}\|(\veps\ptl_{\gx})^{\gamma}f(\gx)\|^2_{L^2(\br^2)}<+\infty,
\end{equation*}
is frequently used, see for instance \cite{arbunich2018rigorous, giannoulis2008interaction, rauch2012hyperbolic}. As shown in Theorem \ref{mainthm}, we do the error estimates in the microscopic scales with the regular Sobolev spaces. Actually, the two treatments are essentially equivalent. Indeed, due to the linearity, the equation remain the same under rescaling $\widetilde{\psi}=\frac{\psi}{\veps}$. The following identity shows the equivalence
\begin{equation*}
	\|\widetilde{\psi}(\gx,T)-\sum_{j=1}^2e^{i\sqrt{E_D}\frac{T}{\veps}}\alpha_j(\gx,T) \Phi_j(\frac{\gx}{\veps})\|_{H^s_{\veps}(\br^2)}=\|\psi(\bx,t)-\sum_{j=1}^2 e^{i\sqrt{E_D}t}\veps \alpha_{j}(\veps\bx,\veps t) \Phi_j(\bx)\|_{H^s(\br^2)}.
\end{equation*}

\emph{Proof of Theorem \ref{mainthm}.}
We follow the standard procedure for wave packet problems. Namely, we first spectrally decompose the error by the Floquet-Bloch theory. Then the spectral components are estimated separately.

Recalling the completeness of Bloch modes of $\ml^A$ in $L^2(\br^2)$, we have
\begin{equation}\label{fb:eta}
	\eta(\bx, t)=\frac{1}{|\Omega^*|}\sum_{b\ge 1}\int_{\Omega^*}\tilde{\eta}_b(\bk, t)\Phi_b(\bx;\bk)d\bk,
\end{equation}
where the error component
\begin{equation}
    \tilde{\eta}_b(\bk, t)=\bla\Phi_b(\bx;\bk),~\eta(\bx, t)\bra.
\end{equation}

Then, for any $b\ge1$, $\tilde{\eta}_b(\bk, t)$ satisfies
\begin{align} \label{etatilde}   \partial^2_{t}\tilde{\eta}_b+E_b(\bk)\tilde{\eta}_b+\bla\Phi_b(\bx;\bk),~ \veps\mlb\eta\bra = e^{i\sqrt{E_D}t}\bla\Phi_b(\bx;\bk),F_1(\bx,t)+F_2(\bx,t)\bra
\end{align}
with initial condition
\begin{equation}
\begin{split}
   \tilde{\eta}_b(\bk,0) &= 0,\\
   \partial_t\tilde{\eta}_b(\bk,t)|_{t=0}   &= \bla\Phi_b(\bx;\bk),F_0(\bx)\bra.
\end{split}
\end{equation}
By Duhamel's principle, we rewrite \eqref{etatilde} as the integral form
\begin{align}\label{Gb}
\nonumber	\tilde{\eta}_b(\bk,t)  =&~~~\frac{e^{i\sqrt{E_b(\bk)}t}-e^{-i\sqrt{E_b(\bk)}t}}{2i\sqrt{E_b(\bk)}}\bla\Phi_b(\bx;\bk),~F_0(\bx)\bra \\
\nonumber    & +\int_{0}^{t}\frac{e^{i\sqrt{E_b(\bk)}(t-\tau)}-e^{-i\sqrt{E_b(\bk)}(t-\tau)}}{2i\sqrt{E_b(\bk)}} e^{i\sqrt{E_D}\tau}\bla\Phi_b(\bx;\bk),~F_1(\bx,\tau)\bra d\tau\\
\nonumber    & +\int_{0}^{t}\frac{e^{i\sqrt{E_b(\bk)}(t-\tau)}-e^{-i\sqrt{E_b(\bk)}(t-\tau)}}{2i\sqrt{E_b(\bk)}} e^{i\sqrt{E_D}\tau}\bla\Phi_b(\bx;\bk),~F_2(\bx,\tau)\bra d\tau\\
\nonumber    & +\int_{0}^{t}\frac{e^{i\sqrt{E_b(\bk)}(t-\tau)}-e^{-i\sqrt{E_b(\bk)}(t-\tau)}}{2i\sqrt{E_b(\bk)}} \bla\Phi_b(\bx;\bk),~\veps\nabla\cdot\big[\nabla\cdot\big(\kappa(\veps \bx)B(\bx)\big)\eta(\bx,\tau)\big]\bra d\tau\\
	:=&~ \tilde{G_0}_b(\bk,t)+\tilde{G_1}_b(\bk,t)+\tilde{G_2}_b(\bk,t)+\tilde{\mathcal{Q}}_b[\eta](\bk,t).
\end{align}
Especially, when $b=1$ and $\bk=\boldsymbol{0}$,
\begin{align}\label{G10}
\nonumber    \tilde{\eta}_1(\boldsymbol{0},t) =& ~~~t\bla\Phi_1(\bx;\boldsymbol{0}),~F_0(\bx)\bra \\
\nonumber    & +\int_{0}^{t}(t-\tau) e^{i\sqrt{E_D}\tau}\bla\Phi_1(\bx;\boldsymbol{0}),~F_1(\bx,\tau)\bra d\tau +\int_{0}^{t}(t-\tau) e^{i\sqrt{E_D}\tau}\bla\Phi_1(\bx;\boldsymbol{0}),~F_2(\bx,\tau)\bra d\tau\\
\nonumber    & +\int_{0}^{t}(t-\tau) \bla\Phi_1(\bx;\boldsymbol{0}),~\veps\nabla\cdot\big[\nabla\cdot\big(\kappa(\veps \bx)B(\bx)\big)\eta(\bx,\tau)\big]\bra d\tau\\
:=&~ \tilde{G_0}_1(\boldsymbol{0},t)+\tilde{G_1}_1(\boldsymbol{0},t)+\tilde{G_2}_1(\boldsymbol{0},t),
\end{align}
the fourth term on the right hand side vanishes since $\Phi_1(\bx; \boldsymbol{0})$ is a constant.

From \eqref{Gb}-\eqref{G10}, we can see that the error components $\tilde{\eta}_b(\bk,t)$ appear in a very different way from those in \cite{fefferman2014wave}. First, $\frac1{\sqrt{E_b(\bk)}}$ appears in the error equation which brings secular terms at $b=1$ and $\bk=\boldsymbol{0}$ as shown in \eqref{G10}. The other main difference is the presence of implicit terms  $\tilde{\mathcal{Q}}_b[\eta](\bk,t)$ which are caused by the effects of modulation/perturbation. Consequently, more efforts and techniques are desired. In our analysis below, we deal with the implicit terms by Gronwall's inequality after proving the boundedness of the operator $\frac{\mlb}{\sqrt{\ml^A}}$ in $H^s(\br^2)$. To handle the singularity, we carefully estimate the error components near and away from the singularity separately.

Note that if $\kappa(\veps\bx)\equiv 0$, which happened similarly in Fefferman and Weinstein's work \cite{fefferman2014wave} on the Schr\"odinger equation, then $\tilde{\mathcal{Q}}_b[\eta](\bk,t),~ b\ge1$ vanished. As a consequence, each $\tilde{\eta}_b(\bk, t)$ is explicitly represented without coupling to other components.

In general, \eqref{Gb} and \eqref{G10} imply that $\eta(\bx, t)$ satisfies the integral equation
\begin{align}\label{eta:bloch}
\nonumber	\eta(\bx, t)=&~\frac{1}{|\Omega^*|}\sum_{b\ge 1}\int_{\Omega^*}\Big[\tilde{G_0}_b(\bk, t)+\tilde{G_1}_b(\bk, t)+\tilde{G_2}_b(\bk, t)+\tilde{\mathcal{Q}}_b[\eta](\bk, t)\Big]\Phi_{b}(\bx;\bk)d\bk \\
	:=&~ G_0(\bx,t)+G_1(\bx,t)+G_2(\bx,t)+\mathcal{Q}[\eta](\bx,t),
\end{align}
and more precisely, $G_j(\bx,t),~j=0,~1,~2$ and $\mathcal{Q}[\eta](\bx,t)$ satisfy
\begin{equation} \label{sequence}
G_j(\bx,t)=\frac{1}{|\Omega^*|}\sum_{b\ge 1}\int_{\Omega^*}\tilde{G_j}_b(\bk,t)\Phi_b(\bx;\bk)d\bk,\quad\mathcal{Q}[\eta](\bx,t)=\frac1{|\Omega^*|}\sum_{b\ge1}\int_{\Omega^*} \tilde{\mathcal{Q}}_b[\eta](\bk,t) \Phi_b(\bx;\bk) d\bk.
\end{equation}
To achieve the error bound of $\eta(\bx, t)$ from \eqref{eta:bloch}, we require the following two important propositions.

\begin{proposition}\label{prop:q}
According to \eqref{Gb}-\eqref{sequence} and the Floquet-Bloch theory \eqref{fb:f}, for any $t\ge0$ and given the integer $s\ge0$, one can get
	\begin{equation}\label{estimate:q}
	\|\mathcal{Q}[\eta](\bx, t)\|^2_{H^s(\br^2)}\le C\veps^2 t \int_{0}^{t}\|\eta(\bx,\tau)\|^2_{H^s(\br^2)} ~d\tau.
	\end{equation}
	
\end{proposition}

\emph{Proof.} We first recall the well-known results on the Riesz transform $\nabla\big(\ml^{A}\big)^{-\frac12}$ which ensures the boundedness of the operator $\nabla\big(\ml^{A}\big)^{-\frac12}$ from $L^2(\br^2)$ to $L^2(\br^2)$, see \cite{auscher2002solution, auscher1998square, hofmann2003p} and the references therein for details. Since the operator $\ml^A$ is self-adjoint, for any two functions $f\in H^1(\br^2)$, $\varphi\in \mathcal{S}(\br^2)$ with $\|\varphi\|_{L^2(\br^2)}=1$, one can obtain the following estimate by Riesz transform in the dual case,
\begin{align}
\nonumber	\bav\bla\big(\ml^{A}\big)^{-\frac12}\nabla f,~\varphi\bra_{L^2(\br^2)}\bav =&~ \bav\bla f,~\nabla{\big(\ml^{A}\big)}^{-\frac12}\varphi\bra_{L^2(\br^2)}\bav \\
\nonumber	\le&~ \|f\|_{L^2(\br^2)}\cdot C\|\varphi\|_{L^2(\br^2)} \\
\le&~ C\|f\|_{L^2(\br^2)}.
\end{align}

Then, for any $t\ge0$, we take the absolute value of $\tilde{\mathcal{Q}}_b[\eta](\bk, t)$,
\begin{equation*}
\begin{split}
\bav\tilde{\mathcal{Q}}_b[\eta](\bk, t)\bav =&~\bav\int_{0}^{t}\frac{\sin\big(\sqrt{E_b(\bk)}(t-\tau)\big)}{\sqrt{E_b(\bk)}} \bla\Phi_b(\bx;\bk),~ \veps\nabla\cdot\big[\nabla\cdot\big(\kappa(\veps \bx)B(\bx)\big)\eta(\bx,\tau)\big]\bra d\tau\bav \\
\le&~\int_{0}^{t} \bav\bla {\big(\ml^A\big)}^{-\frac12}\Phi_b(\bx;\bk),~\veps \nabla\cdot\big[\nabla\cdot\big(\kappa(\veps \bx)B(\bx)\big)\eta(\bx,\tau)\big]\bra \bav d\tau \\
=& ~\int_{0}^{t} \bav\bla \Phi_b(\bx;\bk),~\veps{\big(\ml^{A}\big)}^{-\frac12}\nabla\cdot\big[\nabla\cdot\big(\kappa(\veps \bx)B(\bx)\big)\eta(\bx,\tau)\big]\bra \bav d\tau.
\end{split}
\end{equation*}
By the Floquet-Bloch theory and Minkowski's integral inequality, it follows that
\begin{equation*}
\begin{split}
\|\mathcal{Q}[\eta](\bx, t)\|^2_{L^2(\br^2)} =&~\frac{1}{|\Omega^*|}\sum_{b\ge 1}\int_{\Omega^*}\bav\tilde{\mathcal{Q}}_b[\eta](\bk, t)\bav^2d\bk \\
\le&~\frac{\veps^2}{|\Omega^*|}\sum_{b\ge 1}\Big[\int_0^t \big(\int_{\Omega^*}\bav\bla \Phi_{b}(\bx;\bk),{\big(\ml^{A}\big)}^{-\frac12}\nabla\cdot\big[\nabla\cdot\big(\kappa(\veps \bx)B(\bx)\big)\eta(\bx,\tau)\big]\bra\bav^2 d\bk\big)^{\frac12}d\tau\Big]^2 \\
\le& ~C\veps^2 t\int_{0}^{t} \sum_{b\ge 1}\int_{\Omega^*} \bav\bla \Phi_{b}(\bx;\bk),{\big(\ml^{A}\big)}^{-\frac12}\nabla\cdot\big[\nabla\cdot\big(\kappa(\veps \bx)B(\bx)\big)\eta(\bx,\tau)\big]\bra\bav^2 ~d\bk ~d\tau \\
\le&~C\veps^2 t \int_{0}^{t}\|{\big(\ml^{A}\big)}^{-\frac12}\nabla\cdot\big[\nabla\cdot\big(\kappa(\veps \bx)B(\bx)\big)\eta(\bx,\tau)\big]\|^2_{L^2(\br^2)} ~d\tau\\
\le&~C\veps^2 t \int_{0}^{t}\|\eta(\bx,\tau)\|^2_{L^2(\br^2)} ~d\tau.
\end{split}
\end{equation*}
Further, we obtain the estimate in $H^s(\br^2)$ for any integer $s\ge1$,
\begin{equation*}
\begin{split}
	\|\mathcal{Q}[\eta](\bx, t)\|^2_{H^s(\br^2)} \approx&~\frac{1}{|\Omega^*|}\sum_{b\ge 1}  \int_{\Omega^*}(1+E_b(\bk))^s~ \bav\tilde{\mathcal{Q}}_b[\eta](\bk, t)\bav^2d\bk \\
	\le& ~C \veps^2 t\int_0^t\sum_{b\ge 1}\int_{\Omega^*}  \frac{(1+E_b(\bk))^s}{E_b(\bk)} \bav\bla\Phi_b(\bx;\bk),~ \nabla\cdot\big[\nabla\cdot\big(\kappa(\veps \bx)B(\bx)\big)\eta(\bx,\tau)\big]\bra \bav^2 d\bk~d\tau \\
	\le&~C\veps^2 t \int_{0}^{t} \|\eta(\bx,\tau)\|^2_{L^2(\br^2)}+\|\nabla\cdot\big[\nabla\cdot\big(\kappa(\veps\bx)B(\bx)\big)\eta(\bx,\tau)\big]\|^2_{H^{s-1}(\br^2)} ~d\tau \\
	\le&~C\veps^2 t \int_{0}^{t}\|\eta(\bx,\tau)\|^2_{H^s(\br^2)} ~d\tau.~
\end{split}
\end{equation*}

Now we turn to estimates of the first three terms in \eqref{eta:bloch} which we conclude in the following Proposition.

\begin{proposition}\label{prop:2}
	As $G_j(\bx,t),~j=0,~1,~2$ are defined in \eqref{sequence}, for any $s\ge0$, $\rho>0$, $0<\nu<1$ and $\veps>0$ sufficiently small, the following three statements hold
	\begin{eqnarray}
		\sup_{0\le t\le \rho\veps^{^{-1}}}\|G_0(\bx,t)\|_{H^s(\br^2)} &\leq& C\veps, \label{G0}\\
		\sup_{0\le t\le \rho\veps^{^{-1}}}\|G_1(\bx,t)\|_{H^s(\br^2)} &\leq& C\veps^{1-\nu}, \label{G1}\\
		\sup_{0\le t\le \rho\veps^{^{-1}}}\|G_2(\bx,t)\|_{H^s(\br^2)} &\leq& C\veps, \label{G2}
	\end{eqnarray}
	where each constant $C$ is independent of $\veps$.
\end{proposition}
We shall prove this proposition in the subsequential section.
\newline

Therefore, thanks to the above Proposition \ref{prop:q}-\ref{prop:2}, we can now prove Theorem \ref{mainthm}. For any $s\ge0$, $0\le t\le\rho\veps^{^{-1}}$ and $0<\nu<1$,
\begin{align}
\nonumber \|\eta(\bx, t)\|^2_{H^s(\br^2)} &\le 2\|G_0(\bx, t)+G_1(\bx, t)+G_2(\bx, t)\|^2_{H^s(\br^2)}+2\|\mathcal{Q}[\eta](\bx, t)\|^2_{H^s(\br^2)} \\
\nonumber &\leq C_1\veps^{2-2\nu}+C_2\veps^2 t\int_{0}^{t}\|\eta(\bx, \tau)\|^2_{H^{s}(\br^2)} d\tau.
\end{align}
By utilizing Gronwall's inequality, we obtain
\begin{equation} \label{gronwall:eta}
\|\eta(\bx, t)\|^2_{H^s(\br^2)} \leq C_1\veps^{2-2\nu}\exp(C_2\veps^2 t^2),
\end{equation}
and finally
\begin{equation}
\sup_{0\le t \le \rho\veps^{^{-1}}}\|\eta(\bx, t)\|_{H^s(\br^2)} \leq C\veps^{1-\nu}.
\end{equation}

To complete the proof of  Theorem \ref{mainthm}, we only need to justify \eqref{G0}-\eqref{G2}  of which the detailed verifications will be given successively in the next section.

\section{Proof of Proposition \ref{prop:2}}\label{proof:prop}

In this section, we shall give the detailed proof of Proposition \ref{prop:2}. Hereafter, we just suppress the subscript of $\gk_*$ as $\gk$ for simplicity. For the convenience of our proof, we set $\lambda=\frac{|\gk|}2$ which ensures the lower positive bound in \eqref{eigenvalue:lowbd}, and let $0<q_1<|\gk|-\lambda$ in the Proposition \ref{dirac:modes}. Before proceeding further, we present several results concluded in Proposition \ref{propdirac} which will be frequently used in our proof.
\begin{proposition}\label{propdirac}
	Let $\Gamma(\gx)\in \mathcal{S}(\br^2)$ and $\Psi(\bx)\in C^{\infty}(\br^2)\cap L_{\gk}^2(\br^2/\Lambda)$. Then, the following statements hold
\begin{align}\label{poissonsummation}
	\nonumber & \bla\Phi_b(\bx;\bk),\Gamma(\veps \bx)\Psi(\bx)\big\rangle \\
	=&~ \veps^{-2}\frac{1}{|\Omega|}\int_{\Omega}\overline{\Phi_b(\bx;\bk)}\sum_{\bm\in \mathbb{Z}^2}e^{i(m_1\bk_1+m_2\bk_2+\bk-\mathbf{K})\cdot\bx}\hat{\Gamma}(\frac{m_1\bk_1+m_2\bk_2+\bk-\mathbf{K}}{\veps})\Psi(\bx)d\bx,
\end{align}
and	
\begin{equation}\label{bound:Gamma1}
	\big|\hat{\Gamma}(\boldsymbol{\xi})\big| \leq C\frac{1}{~|\boldsymbol{\xi}|^N}\|\Gamma(\gx)\|_{W^{N,1}(\br^2)}\le C\frac{1}{~|\boldsymbol{\xi}|^N},\quad \forall~\boldsymbol{\xi}\neq \boldsymbol{0}.
\end{equation}
Further, for any $\bk\in\Omega^*$, $\frac12<\nu_1<1$ and $\veps>0$ sufficiently small,
	\begin{eqnarray}
	&|m_1\bk_1+m_2\bk_2+\bk-\gk|\ge C(1+|\bm|),\quad \forall~\bm\in\mathbb{Z}^2,~|\bk-\gk|\geq q_1;& \label{bm:1}\\
	&|m_1\bk_1+m_2\bk_2+\bk-\mathbf{K}|\geq C|\bm|,\quad \forall~\bm\neq(0,0),~|\bk-\gk|<\veps^{\nu_1};& \label{bm:2}\\
	&|m_1\bk_1+m_2\bk_2+\bk-\gk|\ge C\veps^{\nu_1}(1+|\bm|),\quad \forall~\bm\in\mathbb{Z}^2,~|\bk-\gk|\geq \veps^{\nu_1}, \label{bm:3}&
	\end{eqnarray}
where each $C$ is a generic constant.
	
\end{proposition}
The detailed proof is omitted, we refer the readers to \cite{fefferman2014wave} for complete discussions.

Utilizing the Poisson-Summation \eqref{poissonsummation} and the property \eqref{bound:Gamma1} in Proposition \ref{propdirac}, we can obtain the following lemma.
\begin{lemma} \label{lemmaGj1}
Let $\tilde{G_j}_1(\bk,t),~j=0,~1,~2$ be defined in \eqref{Gb}. For any $\bk\in \Omega^*$, $|\bk|<\frac{|\gk|}2$, $N>2$ and $0\le t\le\rho\veps^{^{-1}}$,
	\begin{equation}
	\chi(|\bk|<\frac{|\gk|}2)~|\tilde{G_1}_1(\bk,t)| \le C\veps^{N-2},\quad \chi(|\bk|<\frac{|\gk|}2)~|\tilde{G_j}_1(\bk,t)| \le C\veps^{N-1}\quad j=0,~2.
	\end{equation}
\end{lemma}

\emph{Proof}. Without loss of generality, we only consider the case of $j=1$, the other two cases $j=0,~2$ can be treated similarly. For the convenience, we have to introduce the following notations,
\begin{equation}\label{gammapsi}
\begin{split}
	\Gamma_{11}=-2i\sqrt{E_D}\ptl_T\alpha_1,~\Gamma_{12}=-2i\sqrt{E_D}\ptl_T\alpha_2,~\Gamma_{13}=\nabla_{\gx}&\alpha_1,~\Gamma_{14}=\nabla_{\gx}\alpha_2,~\Gamma_{15}=\kappa\alpha_1,~\Gamma_{16}=\kappa\alpha_2; \\
	\Psi_{11}=\Phi_1,~\Psi_{12}=\Phi_2,~\Psi_{13}=i\mathcal{A}\Phi_1,~\Psi_{14}=i\mathcal{A}\Phi_2&,~\Psi_{15}=-\ml^B \Phi_1,~\Psi_{16}=-\ml^B \Phi_2.
\end{split}
\end{equation}
Thanks to Proposition \ref{schwartz}, we can conclude that
$\Gamma_{1r}=\Gamma_{1r}(\gx,T)\in \mathcal{S}(\br^2)$ for any $T\ge0$, and
$\Psi_{1r}=\Psi_{1r}(\bx) \in C^{\infty}(\br^2)\cap L^2_{\gk}(\br^2/\Lambda),~r\in\{1,\cdots,6\}$.
Then, Poisson-Summation \eqref{poissonsummation} yields that
\begin{align}
\nonumber & \bla\Phi_1(\bx; \bk),\Gamma_{1r}(\veps \bx, \veps \tau)\cdot\Psi_{1r}(\bx)\big\rangle \\
=&~ \veps^{-2}\frac{1}{|\Omega|}\int_{\Omega}\overline{\Phi_1(\bx;\bk)}\sum_{\bm\in \mathbb{Z}^2}e^{i(m_1\bk_1+m_2\bk_2+\bk-\mathbf{K})\cdot\bx}\hat{\Gamma}_{1r} (\frac{m_1\bk_1+m_2\bk_2+\bk-\mathbf{K}}{\veps},\veps \tau)\cdot\Psi_{1r}(\bx)d\bx.
\end{align}

If $|\bk|<\frac{|\gk|}2$, i.e., $|\bk-\gk|\ge q_1$, there exists a constant $C>0$ by \eqref{bm:1} such that
\begin{equation}
	|m_1\bk_1+m_2\bk_2+\bk-\mathbf{K}|\ge C(|\bm|+1), \quad \forall~\bm\in\mathbb{Z}^2.
\end{equation}

According to \eqref{bound:Gamma1},
\begin{eqnarray}
\nonumber  \bav\bla\Phi_{1}(\bx;\bk),\Gamma_{1r}(\veps \bx, \veps \tau)\cdot\Psi_{1r}(\bx)\bra\bav \le  C\veps^{-2}\sum_{\bm\in\mathbb{Z}^2}\frac{\veps^N}{(|\bm|+1)^N}.
\end{eqnarray}

Since $F_1(\bx,\tau)=\veps^2\sum\limits_{r=1}^6\Gamma_{1r}(\veps \bx, \veps \tau)\cdot\Psi_{1r}(\bx)$, we get when $N>2$, $t\in[0,\rho \veps^{^{-1}}]$,
\begin{align}
\nonumber	&\chi(|\bk|<\frac{|\gk|}2)~|\tilde{G_1}_1(\bk,t)| \\
\nonumber	=&~ \chi(|\bk|<\frac{|\gk|}2)~\bav\int_{0}^{t}\frac{\sin\big(\sqrt{E_1(\bk)}(t-\tau)\big)}{\sqrt{E_1(\bk)}} e^{i\sqrt{E_D}\tau}\bla\Phi_1(\bx;\bk),~F_1(\bx,\tau)\bra d\tau \bav\\
\nonumber	\le&~ C\chi(|\bk|<\frac{|\gk|}2)~\int_0^t (t-\tau) \veps^2 \bav\bla\Phi_{1}(\bx;\bk),\Gamma_{1r}(\veps \bx, \veps \tau)\cdot\Psi_{1r}(\bx)\bra\bav d\tau \\
	\le&~ C\veps^{N-2}.
\end{align}

In the following justifications, we will frequently use the consequence in Proposition \ref{schwartz} and \ref{propdirac}.

\subsection{Proof of \eqref{G0}}

Let $N>2$, then by Lemma \ref{lemmaGj1},
\begin{equation}
\chi(|\bk|<\frac{|\gk|}2)~|\tilde{G_0}_1(\bk,t)| \le C\veps.
\end{equation}
\eqref{eigenvalue:lowbd}, \eqref{Gb}-\eqref{sequence} imply
\begin{align}
\nonumber & \chi(|\bk|\ge\dlt_{b,1}\frac{|\gk|}2)~|\tilde{G_0}_b(\bk,t)| \\
\nonumber =&~ \chi(|\bk|\ge\dlt_{b,1}\frac{|\gk|}2)\Big|\frac{e^{i\sqrt{E_b(\bk)}t}-e^{-i\sqrt{E_b(\bk)}t}}{2i\sqrt{E_b(\bk)}}~\veps^2
\bla\Phi_b(\bx;\bk),~\partial_T\alpha_j(\veps \bx,0)\Phi_j(\bx)\bra\Big| \\
\nonumber \leq&~ C\chi(|\bk|\ge\dlt_{b,1}\frac{|\gk|}2)~\veps^2\bav\bla\Phi_b(\bx;\bk),~\partial_T\alpha_j(\veps \bx,0)\Phi_j(\bx)\bra\bav.
\end{align}

Therefore, for any $t\ge0$, we acquire that
\begin{align}\label{P0:1}
\nonumber  & ~\|G_0(\bx,t)\|^2_{H^{s}(\br^2)} \\
\nonumber  \approx&~ \|\frac{1}{|\Omega^*|}\sum_{b\ge1}\int_{\Omega^*}(1+E_b(\bk))^s |\tilde{G_0}_b(\bk,t)|^2 d\bk \\
\nonumber   \le&~ C\int_{\Omega^*}\chi(|\bk|<\frac{|\gk|}2)~|\tilde{G_0}_1(\bk,t)|^2~d\bk+C\sum_{b\geq1}\int_{\Omega^*}(1+E_b(\bk))^s~\chi(|\bk|\ge\dlt_{b,1}\frac{|\gk|}2)~|\tilde{G_0}_b(\bk,t)|^2d\bk \\
\nonumber \leq&~ C\veps^2+C\veps^4\sum_{b\geq1}\int_{\Omega^*}(1+E_b(\bk))^s~\chi(|\bk|\ge\dlt_{b,1}\frac{|\gk|}2)~\bav\bla\Phi_b(\bx;\bk),~\partial_T\alpha_j(\veps \bx,0)\Phi_j(\bx) \bra\bav^2 d\bk \\
\leq&~ C\veps^2 (1+\|\balpha_0(\mathbf{X})\|^2_{H^{s+1}(\br^2)}).
\end{align}
Thus \eqref{G0} is justified.

\subsection{Proof of \eqref{G2}}

As in the case of the proof above, it is straightforward to show that
\begin{align}\label{G2b:1}
\nonumber &~\chi(|\bk|\ge\dlt_{b,1}\frac{|\gk|}2)~|\tilde{G_2}_b(\bk,t)| \\
\nonumber =&~ \chi(|\bk|\ge\dlt_{b,1}\frac{|\gk|}2)~\bav \int_{0}^{t}\frac{e^{i\sqrt{E_b(\bk)})(t-\tau)}-e^{-i\sqrt{E_b(\bk)}(t-\tau)}}{2i\sqrt{E_b(\bk)}}
e^{i\sqrt{E_D}\tau}\bla\Phi_b(\bx;\bk),F_2(\bx,\tau)\bra d\tau\bav \\
	\leq&~ C \chi(|\bk|\ge\dlt_{b,1}\frac{|\gk|}2)\int_{0}^t \bav\bla\Phi_b(\bx;\bk),F_2(\bx,\tau)\bra\bav d\tau.
\end{align}
Moreover, in the calculations below we will use the fact that for any $0\le \tau\le\rho\veps^{^{-1}}$, $s\ge0$,
$$\|F_2(\bx,\tau)\|_{H^{s}(\br^2)}\le C\veps^2.$$

Then, thanks to Lemma \ref{lemmaGj1}, we get for $0\le t \le \rho\veps^{^{-1}}$,
\begin{align}
\nonumber &~\|G_2(\bx, t)\|^2_{H^{s}(\br^2)} \\
\nonumber \leq &~ C\int_{\Omega^*}\chi(|\bk|<\frac{|\gk|}2)~|\tilde{G_2}_1(\bk,t)|^2 d\bk+C\sum_{b\geq1}\int_{\Omega^*}(1+E_b(\bk))^s~\chi(|\bk|\ge\dlt_{b,1}\frac{|\gk|}2)~|\tilde{G_2}_b(\bk,t)|^2d\bk\\
\nonumber \leq&~ C \veps^2+C  t \int_{0}^t \sum_{b\geq1}\int_{\Omega^*}(1+E_b(\bk))^s~\chi(|\bk|\ge\dlt_{b,1}\frac{|\gk|}2)~\bav\bla\Phi_b(\bx;\bk),F_2(\bx,\tau)\bra\bav^2 d\bk~ d\tau \\
\nonumber \leq&~ C \veps^2+C t\int_{0}^t\| F_2(\bx,\tau)\|^2_{H^{s}(\br^2)}d\tau \\
 \leq&~ C \veps^2.
\end{align}
This leads to the result \eqref{G2}, i.e.,
$$\sup_{0\le t\le \rho\veps^{^{-1}}}\|G_2(\bx,t)\|_{H^s(\br^2)} \leq C\veps.$$

\subsection{Proof of \eqref{G1}} \label{esimate:G1}

In this subsection, we turn to the key estimate \eqref{G1}. The main idea of spectral domain decomposition is similar to that presented in \cite{fefferman2014wave} except several considerable modifications due to the differences of the underlying problems. We first divide the Bloch components of $G_1(\bx, t)$ as follows,
\begin{align}
\nonumber   G_{1}(\bx,t) =&~\quad \frac{1}{|\Omega^*|}\int_{\Omega^*}\chi(|\bk|<\frac{|\gk|}2)~\tilde{G_1}_1(\bk, t)\Phi_1(\bx;\bk)d\bk \\
	& +\frac{1}{|\Omega^*|}\sum_{b\ge1}\int_{\Omega^*}\chi(|\bk|\ge\dlt_{b,1}\frac{|\gk|}2)\tilde{G_1}_b(\bk, t)\Phi_b(\bx;\bk)d\bk. \label{G1K}
\end{align}
Since we have assumed that $\veps>0$ sufficiently small satisfying $\veps^{\frac12}<\frac{|\gk|}2$, and $0<q_1<\frac{|\gk|}2$ as mentioned before, the second part on the right hand side of above can be divided into $b\in\{+,-\}$ part ${G_1}_{\mathrm{D}}$ and $b\notin\{+,-\}$ part ${G_1}_{\mathrm{D^C}}$. Further, we decompose ${G_1}_{\mathrm{D}}$ and ${G_1}_{\mathrm{D^C}}$ into their quasi-momentum components near and far away from $\gk$. Specifically, we have the following decomposition,
\begin{align}
\nonumber	{G_1}_{\mathrm{D}}(\bx, t) =&~\quad \frac{1}{|\Omega^*|}\sum_{b\in\{+,-\}}\int_{\Omega^*}\chi(|\bk-\gk|<\veps^{\nu_1})\tilde{G_1}_b(\bk, t)\Phi_b(\bx;\bk)d\bk \\
\nonumber & +\frac{1}{|\Omega^*|}\sum_{b\in \{+,-\}}\int_{\Omega^*}\chi(|\bk-\gk|\geq\veps^{\nu_1},~|\bk|\ge\dlt_{b,1}\frac{|\gk|}2)\tilde{G_1}_b(\bk, t)\Phi_b(\bx;\bk)d\bk \\
:=& ~{G_1}_{\mathrm{D,I}}(\bx, t)+{G_1}_{\mathrm{D,II}}(\bx, t)
\end{align}
and
\begin{align}
\nonumber	{G_1}_{\mathrm{D^C}}(\bx, t) =&~\quad \frac{1}{|\Omega^*|}\sum_{b\notin\{+,-\}}\int_{\Omega^*}\chi(|\bk-\gk|<q_1)\tilde{G_1}_b(\bk, t)\Phi_b(\bx;\bk)d\bk \\
\nonumber   & +\frac{1}{|\Omega^*|}\sum_{b\notin\{+,-\}}\int_{\Omega^*}\chi(|\bk-\gk|\geq q_1,~|\bk|\ge\dlt_{b,1}\frac{|\gk|}2)\tilde{G_1}_b(\bk, t)\Phi_b(\bx;\bk)d\bk \\
:=&~ {G_1}_{\mathrm{D^C,I}}(\bx, t)+{G_1}_{\mathrm{D^C,II}}(\bx, t). \label{divideofG}
\end{align}
Here $\chi$ is the indicator function, $\frac12<\nu_1<1$ is to be determined and $q_1$ is a specified constant independent of $\veps$ such that \eqref{lowerbound} holds.

Then by definition \eqref{l2norm} and \eqref{sobolevnorm}, $G_1(\bx,t)$ can be approximated as follows:
\begin{align}
\nonumber \|G_1(\bx,t)\|^2_{H^{s}(\br^2)} \approx& \quad\int_{\Omega^*}\chi(|\bk|<\frac{|\gk|}2)~|\tilde{G_1}_1(\bk, t)|^2 d\bk \\
\nonumber	& +\|{G_1}_{\mathrm{D,I}}(\bx, t)\|^2_{L^2(\br^2)}+\|{G_1}_{\mathrm{D,II}}(\bx, t)\|^2_{L^2(\br^2)} \\
	& +\|{G_1}_{\mathrm{D^C,I}}(\bx, t)\|^2_{H^{s}(\br^2)}+\|{G_1}_{\mathrm{D^C,II}}(\bx, t)\|^2_{H^{s}(\br^2)}. \label{norm:l2hs}
\end{align}
The estimate of the first part on the right hand side of above directly follows from Lemma \ref{lemmaGj1} when $N\ge3$, and next we verify the error bounds of the other four parts term by term.

Recall that $\tilde{G_1}_b(\bk, t)$ is defined as
\[
	\tilde{G_1}_b(\bk, t)=\int_{0}^{t}\frac{e^{i\sqrt{E_b(\bk)})(t-\tau)}-e^{-i\sqrt{E_b(\bk)}(t-\tau)}}{2i\sqrt{E_b(\bk)}}
	e^{i\sqrt{E_D}\tau} \bla\Phi_b(\bx;\bk),~F_1(\bx,\tau)\bra d\tau,
\]
with
\begin{align}
\nonumber	F_1(\bx,\tau) =& -\veps^2\Big[2i\sqrt{E_D}\ptl_T\alpha_j(\veps\bx,\veps \tau)\Phi_j(\bx)-\nabla_{\gx}\alpha_j(\veps\bx,\veps \tau)\cdot i\mathcal{A}\Phi_j(\bx)+\kappa(\veps\bx)\alpha_j(\veps\bx,\veps \tau)\ml^B\Phi_j(\bx)\Big]\\
	:=& ~\veps^2\sum_{r=1}^{6}\Gamma_{1r}(\veps\bx, \veps \tau)\cdot\Psi_{1r}(\bx),
\end{align}
where $\Gamma_{1r}=\Gamma_{1r}(\gx,T)\in \mathcal{S}(\br^2) ~(\forall~ T\ge0)$, and
$\Psi_{1r}=\Psi_{1r}(\bx) \in C^{\infty}(\br^2)\cap L^2_{\gk}(\br^2/\Lambda),~r\in\{1,\cdots,6\}$ have been defined in \eqref{gammapsi}. By Poisson-Summation stated in \eqref{poissonsummation}, we can directly have
\begin{align}\label{poisson}
\nonumber & \bla\Phi_b(\bx;\bk),\Gamma_{1r}(\veps \bx, \veps \tau)\cdot\Psi_{1r}(\bx)\big\rangle \\
=&~ \veps^{-2}\frac{1}{|\Omega|}\int_{\Omega}\overline{\Phi_b(\bx;\bk)}\sum_{\bm\in \mathbb{Z}^2}e^{i(m_1\bk_1+m_2\bk_2+\bk-\mathbf{K})\cdot\bx}\hat{\Gamma}_{1r} (\frac{m_1\bk_1+m_2\bk_2+\bk-\mathbf{K}}{\veps},\veps \tau)\Psi_{1r}(\bx)d\bx.
\end{align}

Let us start to deal with ${G_1}_{\mathrm{D,I}}(\bx, t)$ at first. For any $\frac12<\nu_1<1$,
\begin{equation*}
	{G_1}_{\mathrm{D,I}}(\bx, t) =  \frac{1}{|\Omega^*|}\sum_{b\in\{+,-\}}\int_{\Omega^*}\chi(|\bk-\gk|<\veps^{\nu_1})\tilde{G_1}_b(\bk, t)\Phi_b(\bx;\bk)d\bk.
\end{equation*}
Then, by the formula \eqref{poisson}, we denote that
\begin{align}
\nonumber &~ \bla\Phi_{\pm}(\bx;\bk),\Gamma_{1r}(\veps \bx, \veps \tau)\cdot\Psi_{1r}(\bx)\bra \\
\nonumber	=& ~\quad\veps^{-2}\frac{1}{|\Omega|}\int_{\Omega}e^{i(\bk-\mathbf{K})\cdot\bx}\hat{\Gamma}_{1r}(\frac{\bk-\mathbf{K}}{\veps},\veps \tau)\cdot\overline{\Phi_{\pm}(\bx;\bk)}\Psi_{1r}(\bx)d\bx \\
\nonumber	&+\veps^{-2}\frac{1}{|\Omega|}\int_{\Omega}\sum_{\substack{\bm\in\mathbb{Z}^2 \\ \bm\neq(0,0)}}e^{i(m_1\bk_1+m_2\bk_2+\bk-\mathbf{K})\cdot\bx}\hat{\Gamma}_{1r}(\frac{m_1\bk_1+m_2\bk_2+\bk-\mathbf{K}}{\veps},\veps \tau)\cdot\overline{\Phi_{\pm}(\bx;\bk)}\Psi_{1r}(\bx)d\bx \\
\nonumber =& ~\quad\veps^{-2}\frac{1}{|\Omega|}\hat{\Gamma}_{1r}(\frac{\bk-\mathbf{K}}{\veps},\veps \tau)\cdot\bla \Phi_{\pm}(\bx;\bk), e^{i(\bk-\mathbf{K})\cdot\bx}\Psi_{1r}(\bx)\bra_{L^2(\Omega)} \\ &+\nonumber \veps^{-2}\frac{1}{|\Omega|}\sum_{\substack{\bm\in\mathbb{Z}^2 \\ \bm\neq(0,0)}}\hat{\Gamma}_{1r}(\frac{m_1\bk_1+m_2\bk_2+\bk-\mathbf{K}}{\veps},\veps \tau) \cdot\bla \Phi_{\pm}(\bx;\bk), e^{i(m_1\bk_1+m_2\bk_2+\bk-\mathbf{K})\cdot\bx}\Psi_{1r}(\bx)\bra_{L^2(\Omega)} \\
	:=&~ \text{I}_1(\bk,\tau)+\text{I}_2(\bk,\tau).
	\label{term2}
\end{align}

According to the statement \eqref{bm:2}, when $\bm\neq(0,0)$, $|\bk-\gk|<\veps^{\nu_1}$ as $\veps>0$ sufficiently small,  there exists a positive number $C>0$ such that
\begin{equation}\label{mk1}
|m_1\bk_1+m_2\bk_2+\bk-\mathbf{K}|\geq C|\bm|.
\end{equation}
Due to the conclusion \eqref{bound:Gamma1}, it follows that for any $0\le \tau  \le\rho\veps^{^{-1}}$,
\begin{eqnarray}
\Big|\chi(|\bk-\gk|<\veps^{\nu_1})\hat{\Gamma}_{1r}(\frac{m_1\bk_1+m_2\bk_2+\bk-\mathbf{K}}{\veps},\veps \tau)\Big|\leq C\frac{\veps^N}{|\bm|^N}, \quad \forall \bm~\neq(0,0).
\end{eqnarray}
Since $\Psi_{1r}(\bx),~r=1,~\cdots,~6$ are smooth and $\{\Phi_b(\bx;\bk)\}_{b\ge1}$ is also complete orthogonal in $L^2(\Omega)$ for any fixed $\bk\in\Omega^*$, we immediately obtain
\begin{eqnarray}\label{I2estimate}
	|\chi(|\bk-\gk|<\veps^{\nu_1})\text{I}_2(\bk,\tau)|\leq C\veps^{N-2}\sum_{\substack{\bm\in\mathbb{Z}^2 \\ \bm\neq(0,0)}}\frac{1}{|\bm|^N},
\end{eqnarray}
and choose $N>2$ to guarantee the above convergence for any $\tau\in[0,~\rho\veps^{^{-1}}]$.

To estimate $\text{I}_1(\bk,\tau)$, we will utilize the expansion of $\Phi_{\pm}(\bx;\bk)$ near $\bk=\gk$ which is given in Proposition \ref{dirac:modes}, i.e., for $|\bk-\gk|<\veps^{\nu_1}\le q_0$,
\begin{eqnarray}\label{Phiasy}
\Phi_{\pm}(\bx;\bk) &=& \frac{e^{i\bkappa\cdot\bx}}{\sqrt{2}}\Big[\frac{\kappa_1+i\kappa_2}{|\bkappa|}\Phi_1(\bx)\pm \Phi_2(\bx)+\mathcal{O}_{H^2_{\gk}(\br^2/\Lambda)}(|\bkappa|)\Big],\  \text{here} \ \bkappa=\bk-\gk.
\end{eqnarray}
Substituting \eqref{Phiasy} into $\text{I}_1(\bk, \tau)$ and by the fact all $\Gamma_{1r}(\cdot, \veps\tau)\in\mathcal{S}(\br^2)$ uniformly for $0\le \tau \le \rho \veps^{^{-1}}$, we directly obtain
\begin{align}\label{I1}
\nonumber   &\sqrt{2}~\hat{\Gamma}_{1r}(\frac{\bk-\mathbf{K}}{\veps},\veps \tau)\cdot\bla \Phi_{\pm}(\bx;\bk),~e^{i(\bk-\gk)\cdot\bx}\Psi_{1r}(\bx)\bra_{L^2(\Omega)} \\
\nonumber   =&~ \frac{\kappa_1+i\kappa_2}{|\bkappa|}\hat{\Gamma}_{1r}(\frac{\bk-\mathbf{K}}{\veps},\veps \tau)\cdot\bla \Phi_1(\bx),~ \Psi_{1r}(\bx)\bra_{L^2(\Omega)} \pm \hat{\Gamma}_{1r}(\frac{\bk-\mathbf{K}}{\veps},\veps \tau)\cdot\bla \Phi_2(\bx), ~\Psi_{1r}(\bx)\bra_{L^2(\Omega)} \\
\nonumber &~ +\mathcal{O}\big(|\bkappa|\big)\\
\nonumber   =&~ \frac{\kappa_1+i\kappa_2}{|\bkappa|}\Big[-2i\sqrt{E_D}\widehat{\partial_T\alpha_1}(\frac{\bkappa}{\veps},\veps \tau)+v_{_F} i\widehat{\partial_{X_1}\alpha_2}(\frac{\bkappa}{\veps},\veps \tau)-v_{_F}\widehat{\partial_{X_2}\alpha_2}(\frac{\bkappa}{\veps},\veps \tau)-\vartheta_{\sharp}\widehat{\kappa\alpha_1}(\frac{\bkappa}{\veps},\veps \tau)\Big] \\
\nonumber & ~\pm \Big[-2i\sqrt{E_D}\widehat{\partial_T\alpha_2}(\frac{\bkappa}{\veps},\veps \tau)+ v_{_F} i\widehat{\partial_{X_1}\alpha_1}(\frac{\bkappa}{\veps},\veps \tau)+v_{_F}\widehat{\partial_{X_2}\alpha_1}(\frac{\bkappa}{\veps},\veps \tau)+\vartheta_{\sharp}\widehat{\kappa\alpha_2}(\frac{\bkappa}{\veps},\veps \tau)\Big] \\
&~ +\mathcal{O}\big(|\bkappa|\big).
\end{align}
In the above calculation, the key is to derive $\bla \Phi_j(\bx),~\Psi_{1r}(\bx)\bra_{L^2(\Omega)},~ j=1,~2, ~r=1,\cdots,6$. By the conclusions stated in Proposition \ref{propinner}, all non-vanishing leading order terms in \eqref{I1} are listed
\begin{align}
\nonumber  -2i\sqrt{E_D}\widehat{\partial_T\alpha_1}(\frac{\bkappa}{\veps},\veps \tau)\bla \Phi_1(\bx), ~\Phi_1(\bx)\bra_{L^2(\Omega)} =&~ -2i\sqrt{E_D}\widehat{\partial_T\alpha_1}(\frac{\bkappa}{\veps},\veps \tau), \\
\nonumber  -2i\sqrt{E_D}\widehat{\partial_T\alpha_2}(\frac{\bkappa}{\veps},\veps \tau)\bla \Phi_2(\bx), ~\Phi_2(\bx)\bra_{L^2(\Omega)} =&~ -2i\sqrt{E_D}\widehat{\partial_T\alpha_2}(\frac{\bkappa}{\veps},\veps \tau), \\
\nonumber  \widehat{\nabla_{\gx}\alpha_2}(\frac{\bkappa}{\veps},\veps \tau)\cdot\bla \Phi_1(\bx),~ i\mathcal{A}\Phi_2(\bx)\bra_{L^2(\Omega)} =&~ v_{_F}\Big( i\widehat{\partial_{X_1}\alpha_2}(\frac{\bkappa}{\veps},\veps \tau)-\widehat{\partial_{X_2}\alpha_2}(\frac{\bkappa}{\veps},\veps \tau)\Big), \\
\nonumber  \widehat{\nabla_{\gx}\alpha_1}(\frac{\bkappa}{\veps},\veps \tau)\cdot\bla \Phi_2(\bx),~ i\mathcal{A}\Phi_1(\bx)\bra_{L^2(\Omega)} =&~ v_{_F}\Big( i\widehat{\partial_{X_1}\alpha_1}(\frac{\bkappa}{\veps},\veps \tau)+\widehat{\partial_{X_2}\alpha_1}(\frac{\bkappa}{\veps},\veps \tau)\Big), \\
\nonumber  -\widehat{\kappa\alpha_1}(\frac{\bkappa}{\veps},\veps \tau)\bla \Phi_1(\bx), ~\ml^B\Phi_1(\bx)\bra_{L^2(\Omega)} =&~ -\vartheta_{\sharp}\widehat{\kappa\alpha_1}(\frac{\bkappa}{\veps},\veps \tau), \\
\nonumber  -\widehat{\kappa\alpha_2}(\frac{\bkappa}{\veps},\veps \tau)\bla \Phi_2(\bx),~ \ml^B\Phi_2(\bx)\bra_{L^2(\Omega)} =&~ \vartheta_{\sharp}\widehat{\kappa\alpha_2}(\frac{\bkappa}{\veps},\veps \tau).
\end{align}

Since $\balpha(\veps \bx, \veps t)$ is the solution of Dirac equation \eqref{dirac}, we have the following identity by Fourier transform,
\begin{equation}
\begin{split}
  &-2i\sqrt{E_D}\widehat{\partial_T\alpha_1}(\frac{\bkappa}{\veps},\veps \tau)+v_{_F} i\widehat{\partial_{X_1}\alpha_2}(\frac{\bkappa}{\veps},\veps \tau)-v_{_F}\widehat{\partial_{X_2}\alpha_2}(\frac{\bkappa}{\veps},\veps \tau)-\vartheta_{\sharp}\widehat{\kappa\alpha_1}(\frac{\bkappa}{\veps},\veps \tau)=0,\\
  & -2i\sqrt{E_D}\widehat{\partial_T\alpha_2}(\frac{\bkappa}{\veps},\veps \tau)+ v_{_F} i\widehat{\partial_{X_1}\alpha_1}(\frac{\bkappa}{\veps},\veps \tau)+v_{_F}\widehat{\partial_{X_2}\alpha_1}(\frac{\bkappa}{\veps},\veps \tau)+\vartheta_{\sharp}\widehat{\kappa\alpha_2}(\frac{\bkappa}{\veps},\veps \tau)=0.
\end{split}
\end{equation}
Then it gives that for any $\tau\in [0,~\rho\veps^{^{-1}}]$ and $|\bk-\gk|<\veps^{\nu_1}$,
\begin{equation}\label{I1estimate}
\Big|\chi(|\bk-\gk|<\veps^{\nu_1})\text{I}_1(\bk,\tau)\Big| \leq C\veps^{-2+\nu_1}.
\end{equation}

Thanks to \eqref{I2estimate} and \eqref{I1estimate},  for all $t\in [0,~\rho\veps^{^{-1}}]$, we have
\begin{align}\label{dlttau11}
\nonumber &~\|{G_1}_{\mathrm{D,I}}(\bx, t)\|_{L^2(\br^2)}^2 \\
\nonumber =&~\frac{1}{|\Omega^*|}\sum_{b\in\{+,-\}}\int_{\Omega^*} \Big|\chi(|\bk-\gk|<\veps^{\nu_1})\tilde{G_1}_b(\bk, t)\Big|^2 d\bk\\
\nonumber =&~ \frac{1}{|\Omega^*|}\sum_{b\in\{+,-\}}\int_{\Omega^*} \Big|\chi(|\bk-\gk|<\veps^{\nu_1}) \\
\nonumber & ~~~\cdot\int_{0}^{t}\frac{e^{i\sqrt{E_b(\bk)})(t-\tau)}-e^{-i\sqrt{E_b(\bk)}(t-\tau)}}{2i\sqrt{E_b(\bk)}}
e^{i\sqrt{E_D}\tau} \bla\Phi_b(\bx;\bk),~F_1(\bx,\tau)\bra d\tau\Big|^2 d\bk \\
\nonumber \le&~ C\sum_{b\in\{+,-\}}\int_{\Omega^*} \Big|\chi(|\bk-\gk|<\veps^{\nu_1})\int_{0}^{t}\veps^2 \big(|\text{I}_1(\bk,\tau)|+|\text{I}_2(\bk, \tau)|\big) d\tau \Big|^2 d\bk \\
\nonumber \le&~ C\sum_{b\in\{+,-\}}\int_{\Omega^*} \Big|\chi(|\bk-\gk|<\veps^{\nu_1})\int_{0}^{t}\veps^2 \veps^{-2+\nu_1}d\tau \Big|^2 d\bk \\
\nonumber \le&~ C\sum_{b\in\{+,-\}}\int_{\Omega^*} \Big|\chi(|\bk-\gk|<\veps^{\nu_1})\veps^{-1+\nu_1} \Big|^2 d\bk \\
    \le&~ C\veps^{-2+4\nu_1}~.
\end{align}
\newline

In the next, we show the estimate of $\|{G_1}_{\mathrm{D,II}}(\bx, t)\|_{L^2(\br^2)}^2$. Recall that
definition
 \begin{equation}
	\|{G_1}_{\mathrm{D,II}}(\bx, t)\|_{L^2(\br^2)}^2=\frac{1}{|\Omega^*|}\sum_{b\in \{+,-\}}\int_{\Omega^*}\chi(|\bk-\gk|\geq\veps^{\nu_1},~|\bk|\ge\dlt_{b,1}\frac{|\gk|}2)|\tilde{G_1}_b(\bk, t)|^2 d\bk.
\end{equation}

Due to the result stated in \eqref{bm:3}, if $|\bk-\gk|\ge \veps^{\nu_1}$, there exists a constant $C$ independent of $\veps$ such that for all $\bm\in\mathbb{Z}^2$,
\[
|m_1\bk_1+m_2\bk_2+\bk-\gk|\ge C\veps^{\nu_1}(1+|\bm|).
\]
Then by invoking \eqref{bound:Gamma1} and Poisson-Summation formula \eqref{poisson}, we can obtain
\begin{eqnarray}\label{dlttau12}
\nonumber && \chi(|\bk-\gk|\ge \veps^{\nu_1})\bav\bla\Phi_{\pm}(\bx;\bk),\Gamma_{1r}(\veps \bx, \veps \tau)\cdot\Psi_{1r}(\bx)\bra\bav \\
	&\leq & C\veps^{-2}\sum_{\bm\in \mathbb{Z}^2}\frac{1}{(1+|\bm|)^N}\veps^{(1-\nu_1)N}. \label{tauneq0}
\end{eqnarray}
Thus, to ensure the convergence of \eqref{dlttau11} and \eqref{dlttau12}, we need to choose $\frac12<\nu_1<1$ i.e., $\nu_1=1-\frac{\nu}2$ which indicates $0<\nu<1$ should not be zero, and let $N$ be large enough such that for any $t\in[0,~\rho\veps^{^{-1}}]$,
\begin{align}
	\|{G_1}_{\mathrm{D,I}}(\bx, t)\|_{L^2(\br^2)} &\leq C \veps^{1-\nu} ,\label{G1di}\\
	\|{G_1}_{\mathrm{D,II}}(\bx, t)\|_{L^2(\br^2)} &\leq C \veps .\label{G1dii}
\end{align}

The above \eqref{G1di} and \eqref{G1dii} imply
\begin{equation}\label{G1D}
\sup_{0\le t \le \rho\veps^{^{-1}}}\|{G_1}_{\mathrm{D}}(\bx, t)\|_{L^2(\br^2)}\le C \veps^{1-\nu}.
\end{equation}
\newline

In order to estimate ${G_1}_{\mathrm{D^C,I}}(\bx, t)$, let us fist introduce that
\begin{equation} \label{G1:DCI}
	\|{G_1}_{\mathrm{D^C,I}}(\bx, t)\|^2_{H^s(\br^2)} \approx \sum_{b\notin\{+,-\}}\int_{\Omega^*}(1+E_b(\bk))^s\chi(|\bk-\gk|<q_1)|\tilde{G_1}_b(\bk, t)|^2 d\bk,
\end{equation}
where $q_1$ is chosen in Proposition \ref{dirac:modes} such that if $|\bk-\gk|< q_1$ and $b\notin\{+,-\}$, the following uniform bound holds
\begin{equation} \label{bd:eb}
	|E_b({\bk})-E_D| \geq C,
\end{equation}
and in addition, $\{E_b(\bk)\}_{b\ge1,~|\bk-\gk|<q_1}$ have a positive lower bound.

By carrying out the integration by part to $\tilde{G_1}_b(\bk, t)$ in $t$, it yields that
\begin{align}
\nonumber  &~\frac{\veps^2}{2i\sqrt{E_b(\bk)}}\int_{0}^{t} e^{i\sqrt{E_b(\bk)}(t-\tau)} e^{i\sqrt{E_D}\tau} \bla\Phi_{b}(\bx;\bk),\Gamma_{1r}(\veps \bx, \veps \tau)\cdot\Psi_{1r}(\bx)\bra d\tau \\
\nonumber  =&~\frac{\veps^2}{2i\sqrt{E_b(\bk)}}e^{i\sqrt{E_b(\bk)}t}\int_{0}^{t} e^{i(\sqrt{E_D}-\sqrt{E_b(\bk)})\tau} \bla\Phi_{b}(\bx;\bk),\Gamma_{1r}(\veps \bx, \veps \tau)\cdot\Psi_{1r}(\bx)\bra d\tau \\
\nonumber    =&~ ~~~\frac{\veps^2}{2\sqrt{E_b(\bk)}}\frac{e^{i\sqrt{E_b(\bk)}t}}{\sqrt{E_D}-\sqrt{E_b(\bk)}}\bla\Phi_{b}(\bx;\bk),\Gamma_{1r}(\veps \bx, 0)\cdot\Psi_{1r}(\bx)\bra  \\
\nonumber    &~ -\frac{\veps^2}{2\sqrt{E_b(\bk)}}\frac{e^{i\sqrt{E_D}t}}{\sqrt{E_D}-\sqrt{E_b(\bk)}}\bla\Phi_{b}(\bx;\bk),\Gamma_{1r}(\veps \bx, \veps t)\cdot\Psi_{1r}(\bx)\bra \\
\nonumber    &~ +\frac{\veps^3}{2\sqrt{E_b(\bk)}}\frac{e^{i\sqrt{E_b(\bk)}t}}{\sqrt{E_D}-\sqrt{E_b(\bk)}}\int_{0}^{t}e^{i(\sqrt{E_D}-\sqrt{E_b(\bk)})\tau}\bla\Phi_{b}(\bx;\bk), \ptl_T \Gamma_{1r}(\veps \bx, \veps \tau)\cdot\Psi_{1r}(\bx)\bra d\tau,
\end{align}
and similarly,
\begin{align}
\nonumber  &~\frac{\veps^2}{2i\sqrt{E_b(\bk)}}\int_{0}^{t}e^{-i\sqrt{E_b(\bk)}(t-\tau)}e^{i\sqrt{E_D}\tau}\bla\Phi_{b}(\bx;\bk),\Gamma_{1r}(\veps \bx, \veps \tau)\cdot\Psi_{1r}(\bx)\bra d\tau \\
%\nonumber  =&~\frac{\veps^2}{2i\sqrt{E_b(\bk)}}e^{-i\sqrt{E_b(\bk)}t}\int_{0}^{t}e^{i(\sqrt{E_b(\bk)}+\sqrt{E_D})\tau}\bla\Phi_{b}(\bx;\bk),\Gamma_{1r}(\veps \bx, \veps \tau)\cdot\Psi_{1r}(\bx)\bra d\tau \\
\nonumber    =&~ ~~~\frac{\veps^2}{2\sqrt{E_b(\bk)}}\frac{e^{-i\sqrt{E_b(\bk)}t}}{\sqrt{E_D}+\sqrt{E_b(\bk)}}\bla\Phi_{b}(\bx;\bk),\Gamma_{1r}(\veps \bx, 0)\cdot\Psi_{1r}(\bx)\bra  \\
\nonumber    &~ -\frac{\veps^2}{2\sqrt{E_b(\bk)}}\frac{e^{i\sqrt{E_D}t}}{\sqrt{E_D}+\sqrt{E_b(\bk)}}\bla\Phi_{b}(\bx;\bk),\Gamma_{1r}(\veps \bx, \veps t)\cdot\Psi_{1r}(\bx)\bra \\
\nonumber    &~ +\frac{\veps^3}{2\sqrt{E_b(\bk)}}\frac{e^{-i\sqrt{E_b(\bk)}t}}{\sqrt{E_D}+\sqrt{E_b(\bk)}}\int_{0}^{t}e^{i(\sqrt{E_D}+\sqrt{E_b(\bk)})\tau}\bla\Phi_{b}(\bx;\bk), \ptl_T \Gamma_{1r}(\veps \bx, \veps \tau)\cdot\Psi_{1r}(\bx)\bra d\tau.
\end{align}

Utilizing \eqref{bd:eb} and the unform lower bound of $E_b(\bk)$ when $|\bk-\gk|<q_1$, we immediately conclude the following estimate holds
\begin{align}
\nonumber  |\tilde{G_1}_b(\bk, t)| \leq& ~~~ C\frac{\veps^2}{\sqrt{E_b(\bk)}}\big(\frac{1}{\sqrt{E_D}+\sqrt{E_b(\bk)}}+\Big|\frac{1}{\sqrt{E_D}-\sqrt{E_b(\bk)}}\Big|\big)\\
\nonumber  & \cdot\big(\Big|\bla\Phi_{b}(\bx;\bk),\Gamma_{1r}(\veps \bx, 0)\cdot\Psi_{1r}(\bx)\bra\Big|+\Big|\bla\Phi_{b}(\bx;\bk),\Gamma_{1r}(\veps \bx, \veps t)\cdot\Psi_{1r}(\bx)\bra\Big|\big) \\
\nonumber &+C\frac{\veps^3}{\sqrt{E_b(\bk)}}\big(\frac{1}{\sqrt{E_D}+\sqrt{E_b(\bk)}}+\Big|\frac{1}{\sqrt{E_D}-\sqrt{E_b(\bk)}}\Big|\big) \\
\nonumber & \cdot\int_{0}^{t}\Big|\bla\Phi_{b}(\bx;\bk), \ptl_T \Gamma_{1r}(\veps \bx, \veps \tau)\cdot\Psi_{1r}(\bx)\bra\Big| d\tau\\
\nonumber  \leq & ~~~C \veps^2 \big(\Big|\bla\Phi_{b}(\bx;\bk),\Gamma_{1r}(\veps \bx, 0)\cdot\Psi_{1r}(\bx)\bra\Big|+\Big|\bla\Phi_{b}(\bx;\bk),\Gamma_{1r}(\veps \bx, \veps t)\cdot\Psi_{1r}(\bx)\bra\Big|\big) \\
& +C \veps^3\int_{0}^{t}\Big|\bla\Phi_{b}(\bx;\bk), \ptl_T \Gamma_{1r}(\veps \bx, \veps \tau)\cdot\Psi_{1r}(\bx)\bra\Big| d\tau.
\end{align}

Then by Minkowski's integral inequality, for any $t\in [0,~\rho\veps^{^{-1}}]$ we have
\begin{align}\label{G1DC1}
\nonumber   &~\|{G_1}_{\mathrm{D^C,I}}(\bx, t)\|_{H^s(\br^2)}^2 \\
\nonumber \approx&~ \sum_{b\notin\{+,-\}} \int_{\Omega^*} (1+E_b(\bk))^s \chi(|\bk-\gk|<q_1)|\tilde{G_1}_b(\bk,t)|^2 d\bk \\
\nonumber   \leq&~ C\veps^4\sum_{b\notin\{+,-\}} \int_{\Omega^*}(1+E_b(\bk))^{s}\Big(\Big|\bla\Phi_{b}(\bx;\bk),\Gamma_{1r}(\veps \bx, 0)\cdot\Psi_{1r}(\bx)\bra\Big|^2+\Big|\bla\Phi_{b}(\bx;\bk),\Gamma_{1r}(\veps \bx, \veps t)\cdot\Psi_{1r}(\bx)\bra\Big|^2 \Big) d\bk \\
\nonumber & +C\veps^6\sum_{b\notin\{+,-\}} \int_{\Omega^*}(1+E_b(\bk))^{s}\Big(\int_{0}^{t}\bav\bla\Phi_b(\bx;\bk),\ptl_T \Gamma_{1r}(\veps \bx, \veps \tau)\cdot\Psi_{1r}(\bx)\bra\bav d\tau\Big)^2d\bk\\
\nonumber \le &~ C\veps^2\|\balpha_0(\gx)\|_{H^{s+1}(\br^2)}^2 +C\veps^6 t \int_0^t \sum_{b\notin\{+,-\}}\int_{\Omega^*}(1+E_b(\bk))^{s}\bav\bla\Phi_b(\bx;\bk),\ptl_T \Gamma_{1r}(\veps \bx, \veps \tau)\cdot\Psi_{1r}(\bx)\bra\bav^2 d\bk ~d\tau \\
\nonumber \leq&~ C\veps^2+ C\veps^6 t^2 \sup_{0\leq \tau\leq \rho\veps^{^{-1}}}\|\ptl_T \Gamma_{1r}(\veps \bx, \veps \tau)\cdot\Psi_{1r}(\bx)\|_{H^{s}(\br^2)}^2 \\
\leq&~ C\veps^{2}.
\end{align}

~\newline

We need a new approach to estimate ${G_1}_{\mathrm{D^C,II}}(\bx, t)$. Indeed, the technique of integration by part in the previous estimate will be invalid, since the spectral band $E_b(\bk)$ varies unclearly when $|\bk-\gk| \geq q_1$ and $|\bk|\ge\dlt_{b,1}\frac{|\gk|}2$ in $\Omega^*$. Here the strategy is similar to that of verifying ${G_1}_{\mathrm{D,II}}$. Main difference is that we need to control under the $H^s$ norm, thus a sharper estimate is carried out by applying $(I+\ml^A)$ on $\Phi_b(\bx;\bk)$ as many times as required. Namely,
for any integer $M\ge0$,
\begin{align}
\nonumber  \bla\Phi_b(\bx;\bk), \Gamma_{1r}(\veps \bx, \veps \tau)\cdot\Psi_{1r}(\bx)\bra &= \frac1{(1+E_b(\bk))^{M}}\bla (I+\mathcal{L}^A)^{M}\Phi_b(\bx;\bk), \Gamma_{1r}(\veps \bx, \veps \tau)\cdot\Psi_{1r}(\bx)\bra \\
\nonumber  &= \frac1{(1+E_b(\bk))^{M}}\bla \Phi_b(\bx;\bk), (I+\mathcal{L}^A)^{M}\big(\Gamma_{1r}(\veps \bx, \veps \tau)\cdot\Psi_{1r}(\bx)\big)\bra.
\end{align}

Note that $(I+\mathcal{L}^A)^{M}\big(\Gamma_{1r}(\veps \bx, \veps \tau)\cdot\Psi_{1r}(\bx)\big)$ is a summation of terms in the form $\Upsilon(\veps \bx, \veps \tau) \Theta(\bx)$ where $\Upsilon(\cdot, \veps \tau) \in \mathcal{S}(\br^2)$ includes derivatives of $\Gamma_{1r}(\veps \bx, \veps \tau)$, and $\Theta(\bx)\in C^\infty(\br^2) \cap L^2_\gk(\br^2/\Lambda)$ contains derivatives of $A(\bx)$ and $\Psi_{1r}(\bx)$. According to \eqref{bm:1} in Proposition \ref{propdirac}, if $\bk\in\Omega^*$ and $|\bk-\gk|\geq q_1$, there exists a constant $C>0$ such that
\begin{equation*}
	|m_1\bk_1+m_2\bk_2+\bk-\gk|\ge C(1+|\bm|),\quad \forall~\bm\in\mathbb{Z}^2.
\end{equation*}

Therefore, by applying \eqref{poissonsummation} and \eqref{bound:Gamma1}, we can conclude that for any $\bk\in\Omega^*,~|\bk-\gk|\geq q_1$ and $\tau \in [0,~\rho\veps^{^{-1}}]$,
\begin{equation}
	 \Big|\bla\Phi_b(\bx;\bk), (I+\mathcal{L}^A)^{M} \Gamma_{1r}(\veps \bx, \veps \tau)\cdot\Psi_{1r}(\bx)\bra\Big| \le C\veps^{N-2}\sum_{\bm\in\mathbb{Z}^2}\frac{1}{(1+|\bm|)^N}~.
\end{equation}
Thus, one can figure out for any $t\in[0,~\rho\veps^{^{-1}}]$
\begin{align*}
\nonumber   &~\chi(|\bk-\gk|\ge q_1,~|\bk|\ge\dlt_{b,1}\frac{|\gk|}2)~ |\tilde{G_1}_b(\bk, t)| \\
   \leq&~ C\chi(|\bk-\gk|\ge q_1,~|\bk|\ge\dlt_{b,1}\frac{|\gk|}2)\veps^2 \frac{\rho\veps^{^{-1}}}{(1+E_b(\bk))^M}\veps^{N-2}\sum_{\bm\in\mathbb{Z}^2}\frac{1}{(1+|\bm|)^N}~.
\end{align*}
Here we require $N>2$ to ensure the convergence of the double summation. This leads to the last estimation as follows
\begin{align}\label{G1DC2}
\nonumber   \|{G_1}_{\mathrm{D^C,II}}(\bx, t)\|_{H^s(\br^2)}^2 \approx&~ \sum_{b\geq1}(1+b)^s \int_{\Omega^*} \chi(|\bk-\gk|\ge q_1,~|\bk|\ge\dlt_{b,1}\frac{|\gk|}2)~|\tilde{G_1}_b(\bk, t)|^2 d\bk \\
\nonumber    \leq&~ C\sum_{b\geq1}(1+b)^{s-2M}\veps^2 \\
\leq&~ C\veps^2,
\end{align}
where we choose  $M>(s+1)/2$.

Finally, \eqref{G1D}, \eqref{G1DC1}, \eqref{G1DC2} and Lemma \ref{lemmaGj1} imply that for any $s\ge0$, $\rho>0$, $0<\nu<1$ and $\veps>0$ sufficiently small
\begin{equation}
	\sup_{0\le t \le \rho\veps^{^{-1}}}\|{G_1}(\bx, t)\|_{H^s(\br^2)}\le C \veps^{1-\nu}.
\end{equation}
This completes  the proof of Proposition \ref{prop:2} and therewith Theorem \ref{mainthm}.

\hfill$\square$

~\\

\noindent{\bf Acknowledgments:}
The authors would thank Prof. Michael I. Weinstein for proposing this interesting problem and for his useful suggestions. This work was supported by National Natural Science Foundation of China NSFC grants $\#11871299$.

\begin{appendices}

\section{Appendix: Wave packets with data spectrally localized near Dirac points}
In this appendix, we consider the case of $\kappa(\veps \bx)\equiv 0$ and aim to obtain an asymptotic solution that would be valid over the time scale up to the order $\mathcal{O}(\frac{1}{\veps^{^{2^-}}})$, which is parallel to the main result in \cite{fefferman2014wave}. To put the analysis in the context,
consider the following wave packet equations with data spectrally localized near Dirac points
\begin{equation} \label{weq:unperturbed}
	\partial_{t}^2\psi-\nabla\cdot A(\bx)\nabla \psi = 0,
\end{equation}
with the same initial conditions \eqref{ini}.

We obtain a similar asymptotic solution in the form of
\begin{equation} \label{solu:unperturb}
	\psi(\bx, t)=e^{i\sqrt{E_D}t}\veps\Big[\alpha_{1}(\veps \bx, \veps t)\Phi_1(\bx)+\alpha_{2}(\veps \bx, \veps t)\Phi_2(\bx)\Big]+\eta(\bx, t).
\end{equation}
Here the envelopes $\alpha_j(\veps \bx, \veps t),~j=1,~2$ satisfy the massless Dirac equation
\begin{equation}\label{dirac:unpertrubed}
\left\{
\begin{aligned}
&i\ptl_{T}\alpha_1-\frac{v_{_F}}{2\sqrt{E_D}}\big(i\ptl_{X_1}-\ptl_{X_2}\big)\alpha_2=0 \\
&i\ptl_{T}\alpha_2-\frac{v_{_F}}{2\sqrt{E_D}}\big(i\ptl_{X_1}+\ptl_{X_2}\big)\alpha_1=0
\end{aligned}~,
\right.
\end{equation}
with initial condition $\alpha_{1}(\gx,0)=\alpha_{10}(\gx),~\alpha_{2}(\gx,0)=\alpha_{20}(\gx)$.

To conclude, we can obtain the following result,
\begin{theorem}\label{mainthmunperturb}
Supposed that $A(\bx)$ is a honeycomb structured material weight defined in Definition \ref{hcdef}, $\Phi_1(\bx)$, $\Phi_2(\bx)$ are the eigenfunctions associated with the Dirac point $( \gk, E_D)$ given in Definition \ref{Dirac pts},  $\alpha_{10}(\gx)$, $\alpha_{20}(\gx)$ are Schwartz functions of $\gx\in \br^2$. Then problem \eqref{weq:unperturbed} \eqref{ini} has a unique solution of the form \eqref{solu:unperturb}, where $\alpha_j(\gx,T),~j=1,~2$ are the solution to the massless system \eqref{dirac:unpertrubed}, and for any $s\ge 0$, $\rho>0$, $0<\nu<1$, a sufficiently small positive parameter $\veps$,
\begin{equation}
\sup_{0\le t \le \rho\veps^{^{-2+\nu}}}\|\eta(\bx, t)\|_{H^s(\br^2)}\le C\veps^{\frac{\nu}{2}}.
\end{equation}
Here $C$ does not depend on $\veps$.
\end{theorem}

Firstly, the well-posedness of the massless Dirac equation \eqref{dirac:unpertrubed} for $T\in [0,\infty)$ can be obtained by the Fourier transform, see \cite{fefferman2014wave}. Since the error estimate is quite similar to that in the proof of Theorem \ref{mainthm}, we omit the detailed repeated calculations by just sketching out the main idea. After substituting \eqref{solu:unperturb} into equation \eqref{weq:unperturbed} and decomposing $\eta(\bx, t)$ into its Floquet-Bloch components, we obtain a similar equation \eqref{Gb} \eqref{eta:bloch} without the implicit term $\tilde{\mathcal{Q}}_b[\eta]$. By carefully dealing with \eqref{G0}-\eqref{G2}, we can improve the estimates as follows
\begin{eqnarray*}
\sup_{0\le t\le \rho\veps^{^{-2+\nu}}}\|G_0(\bx,t)\|_{H^s(\br^2)} &\leq& C\veps, \\
\sup_{0\le t\le \rho\veps^{^{-2+\nu}}}\|G_1(\bx,t)\|_{H^s(\br^2)} &\leq& C\veps^{\frac{\nu}2}, \\
\sup_{0\le t\le \rho\veps^{^{-2+\nu}}}\|G_2(\bx,t)\|_{H^s(\br^2)} &\leq& C\veps^{\nu}.
\end{eqnarray*}
Thus Theorem \ref{mainthmunperturb} is asserted.

\section[Proof of Proposition \ref{schwartz}]{Appendix: Proof of Proposition \ref{schwartz}}

The proof of Proposition \ref{schwartz} is given in this appendix. We first list the standard result on the well-posedness of linear symmetric hyperbolic system, see for instance \cite{kato1975cauchy, racke1992lectures}.

Let $U=U(\bx, t)$ be the unknown $n$-vector-valued function from $(\bx,t)\in\br^d\times[0,+\infty)$ to $\br^n$, and $A_1(\bx,t),~ \cdots,~A_d(\bx,t),~B(\bx,t)$ are $n\times n$-matrix-valued functions. Consider the first-order partial differential equations of the form
\begin{equation}\label{hyperbolic}
	\ptl_{t}U+\sum_{j=1}^{d}A_j(\bx,t)\ptl_{x_j}U+B(\bx,t)U = 0,
\end{equation}
with initial condition $U(\bx,0) = U_0(\bx)$.
The global existence of the solution is shown in the following Proposition.
\begin{proposition} \label{thm:hyper}
	Let $s\in \mathbb{N}$ and $s>\frac{d}{2}+1$. Assume that $A_1, \cdots, A_d\in C_b^{s+1}[~or~C_b^{\infty}]$ are Hermitian matrices, $B\in C_b^{s+1}[~or~C_b^{\infty}]$ and $U_0\in H^s [~or~\cap C^{\infty}]$. Then there exists a unique solution $U\in C^1\big(\br^d\times [0, +\infty)\big)[~or~C^{\infty}\big(\br^d\times[0, +\infty) \big)]$, and for any $0<\rho<\infty$,
	\begin{equation}\label{hsestimate}
	U(\bx,t)\in C^{0}\big([0, \rho], H^s(\br^d)\big)\cap C^{1}\big([0, \rho], H^{s-1}(\br^d)\big).
	\end{equation}
	Moreover, there exist $C_1,~C_2>0$ such that
	\begin{equation}
	\|U(\cdot,t)\|_{H^s(\br^d)}\leq C_1 e^{C_2 t}\|U_0(\bx)\|_{H^s(\br^d)}.
	\end{equation}	
\end{proposition}

To prove Proposition \ref{schwartz}, we rewrite the Dirac equation \eqref{dirac} in the compact form
\begin{equation}
\ptl_T\balpha+A_1\ptl_{X_1}\balpha+A_2\ptl_{X_2}\balpha+B(\gx)\balpha=0,
\end{equation}
where the coefficients $A_1=-\frac{v_{_F}}{2\sqrt{E_D}}\sigma_1,~ A_2=\frac{v_{_F}}{2\sqrt{E_D}}\sigma_2$ are Hermitian matrices and $B(\gx)=-i\frac{\vartheta\kappa(\gx)}{2\sqrt{E_D}}\sigma_3$ is a smooth bounded matrix-valued function.

Note that the initial value $\balpha_0(\gx)=\big(\alpha_{10}(\gx),\alpha_{20}(\gx)\big)^T$ is in Schwartz space. The first conclusion \eqref{wellpose} of Proposition \ref{schwartz} is just a direct consequence of Proposition \ref{thm:hyper}. However, it requires a more delicate estimate to derive the second result \eqref{schwartzestimate} in Proposition \ref{schwartz}.

For any integer $N\ge0$ and the multi-indices $\mathbf{n}\in\mathbb{N}^2$ satisfying $|\mathbf{n}|\le N$, define
\begin{equation}
V_{N}:=(\balpha,~\ptl_{X_1}\balpha,~\ptl_{X_2}\balpha,~\cdots,~\ptl^{\mathbf{n}}_{\gx}\balpha,~\cdots)^T.
\end{equation}
According to the Dirac equation \eqref{dirac}, one can deduce $V_{N}$ satisfies
\begin{equation}
\ptl_{T}V_N+\mathscr{A}_1\ptl_{X_1}V_N+\mathscr{A}_2\ptl_{X_2}V_N+\mathscr{B}(\gx)V_N=0,
\end{equation}
with the initial value
\begin{eqnarray}
V_N(\gx,~0)=(\balpha_0,~\ptl_{X_1}\balpha_0,~\ptl_{X_2}\balpha_0,~\cdots,~\ptl^{\mathbf{n}}_{\gx}\balpha_0,~\cdots)^T,
\end{eqnarray}
where
\begin{eqnarray}
\mathscr{A}_1=
\begin{pmatrix}
A_1 & & & \\
& A_1 & & \\
& & \ddots & \\
& & & A_1
\end{pmatrix},\quad
\mathscr{A}_2=
\begin{pmatrix}
A_2 & & &\\
& A_2 & &\\
& &\ddots &\\
& & & A_2
\end{pmatrix}
\end{eqnarray}
are both block diagonal matrices, and $\mathscr{B}(\gx)$ is a lower triangular matrix with each block entry being a linear combination of $\ptl^{\bm}_{\gx}B(\gx), ~\bm\in \mathbb{N}^2,~|\bm|\leq N$, and thus smooth and bounded.

Therefore, for any $M\in \mathbb{N}$, $(1+|\gx|^2)^{\frac M2}V_N$ satisfies the following system
\begin{equation*}
	\ptl_{T}[(1+|\gx|^2)^{\frac M2}V_N]+\mathscr{A}_1\ptl_{X_1}[(1+|\gx|^2)^{\frac M2}V_N]+\mathscr{A}_2\ptl_{X_2}[(1+|\gx|^2)^{\frac M2}V_N]+\widetilde{\mathscr{B}}(\gx)[(1+|\gx|^2)^{\frac M2}V_N]=0
\end{equation*}
with initial values
\begin{equation}
	(1+|\gx|^2)^{\frac M2}V_N(\gx,~0)=(1+|\gx|^2)^{\frac M2} (\balpha_0,~\ptl_{X_1}\balpha_0,~\ptl_{X_2}\balpha_0,~\cdots,~\ptl^{\mathbf{n}}_{\gx}\balpha_0,~\cdots)^T,
\end{equation}
where $\widetilde{\mathscr{B}}(\gx)=\mathscr{B}(\gx)-\mathscr{A}_1\frac{M X_1}{1+|\gx|^2}-\mathscr{A}_2\frac{M X_2}{1+|\gx|^2}$ is also smooth and bounded.

By Proposition \ref{thm:hyper}, we can conclude for any $s>2,~0<\rho<\infty$,
\begin{equation}
(1+|\gx|^2)^{\frac M2}V_N(\gx, T)\in C\big([0,~\rho], H^s(\br^2)\big)\cap C^1\big([0,~\rho], H^{s-1}(\br^2)\big),
\end{equation}
and further for any $T\ge0$,
\begin{equation}\label{VNestimate}
\|(1+|\gx|^2)^{\frac M2}V_N(\gx, T)\|_{H^s(\br^2)}\leq C_1 e^{C_2T} \|(1+|\gx|^2)^{\frac M2}V_N(\gx, 0)\|_{H^s(\br^2)}.
\end{equation}

Then \eqref{VNestimate} implies $(1+|\gx|^2)^{\frac M2}V_N$ is in H$\ddot{o}$lder space by Sobolev embedding theorem, and for any $\rho>0$,
\begin{equation}
\sup_{T\in[0,\rho],~\gx\in \br^2}\bav (1+|\gx|^2)^{\frac M2}V_N(\gx, T)\bav\leq C.
\end{equation}

If $l\in\mathbb{N}$, one can observe by induction, $\ptl_T^l\balpha$ is a linear combination of $\ptl_{\gx}^{\bm}\balpha,~0\le|\bm|\le l$, with the smooth bounded coefficients. Specifically, $\ptl_{T}^2 \alpha_1,~\ptl_{T}^2 \alpha_2$ satisfy an uncoupled form:
\begin{equation*}
\ptl_{T}^2 \alpha_j-\frac{v_{_F}^2}{4E_D}(\ptl^2_{X_1}+\ptl^2_{X_2})\alpha_j+\frac{\vartheta_{\sharp}^2\kappa(\gx)^2}{4E_D}\alpha_j=0\quad(j=1,~2).
\end{equation*}
In conclusion, one can get for any $l,~M\in\mathbb{N},~\bn\in\mathbb{N}^2$, there exists a positive constant $C$ such that
\begin{equation}
\sup_{T\in[0,\rho],~\gx\in \br^2}\bav(1+|\gx|^2)^{\frac M2}\ptl_{\gx}^{\mathbf{n}}\ptl_T^l \balpha(\gx, T)\bav < \infty.
\end{equation}
This completes the proof.

\hfill$\square$

\end{appendices}

\end{document}